\documentclass[10pt, a4paper]{amsart}
\input xy   
\xyoption{all}
\numberwithin{equation}{subsection}

\theoremstyle{plain}

  \begingroup
        \newtheorem{theorem}[equation]{Theorem}
        \newtheorem{lemma}[equation]{Lemma}
        \newtheorem{proposition}[equation]{Proposition}
        \newtheorem{corollary}[equation]{Corollary}

	\newtheorem{definition}[equation]{Definition}
        
  \endgroup
\theoremstyle{definition}
  \begingroup
        \newtheorem{example}[equation]{Example}
        \newtheorem{comment}[equation]{Comment}
  \endgroup

\newcommand{\diagrama}{\xymatrix}
\newcommand{\mr}[1]{\buildrel {#1} \over \longrightarrow}
\newcommand{\ml}[1]{\buildrel {#1} \over \longleftarrow}
\newcommand{\cc}{\mathcal}

\begin{document}
This is a revised version of \textbf{arXiv.org/math.CT/02008222} 

to appear in JPAA. 

\title{On the representation theory of Galois and Atomic topoi}

\author{Eduardo J. Dubuc}

\maketitle
{\sc introduction} \indent \vspace{1ex}

The notion of a (pointed) \emph{Galois pretopos} (``cat\'egorie
Galoisienne'') was considered originally by 
Grothendieck in \cite{G1} in connection with the fundamental 
group of a scheme. In that paper \emph{Galois theory is conceived as the
axiomatic characterization of the classifying pretopos of a profinite
group $G$}. The fundamental theorem takes the form of a
representation theorem for Galois pretopos  (see \cite{DC} for the explicit
interpretation of this work in terms of filtered unions of categories 
 - the link to filtered inverse limits of
topoi -  and its relation to classical Galois's galois theory).
An important motivation was
pragmatical. The fundamental theorem is tailored to be applied
to the category of
etal coverings of a connected locally noetherian scheme pointed 
with a geometric point over an algebraically closed field. We quote:
``Cette \'equivalence permet donc de 
interpr\'eter les op\'erations courantes sur des rev\^etements en
terms des op\'erations analogues dans $\cc{B}G$, i.e. en terms des
op\'erations \'evidentes sur des ensembles finis o\`u $G$ op\'ere''.    
Later, in collaboration with Verdier (\cite{G2} Ex IV), he considers  the
general notion of pointed Galois Topos  in a series
of commented exercises (specially Ex IV, 2.7.5). There, specific
guidelines are given to develop the theory of classifying topoi of
\emph{progroups}. It is stated therein that 
Galois topoi correspond exactly, as categories, to the full subcategories
generated by \emph{locally constant objects} in connected locally connected
topoi (this amounts to the construction of Galois closures), and that
they classify progroups.  In \cite{M2}, Moerdiejk
developed this program under the light of the \emph{localic group}
concept. He proves the fundamental
theorem (in a rather sketchy way, theorem 3.2 loc.cit.) in the form of
a characterization of 
pointed Galois topoi as the  
classifying topoi of \emph{prodiscrete localic groups}.

In appendix-section \ref{locallyconstantsection} we develop the
theory of 
locally constant objects as defined in \cite{G2} Ex. IX. For the
notion of Galois Topos discussed here see definition 
\ref{galoistopos}. We take from
\cite{B} the idea of presenting the topos of objects split by a cover as a
push-out topos. We show how the existence of Galois closures follows
automatically by the fact that this topos has essential points.  

Connected groupoids are
considered already in \cite{G1} because of the lack of a canonical
point. The groupoid whose objects are all the points and with arrows the
natural transformations, imposes itself
as the natural mathematical object to be considered (although all the
information is already in any one of its vertex groups). The theory is
developed with groups for the sake of simplicity, but the appropriate
formulation of the groupoid version is not straightforward 
\mbox{(see \cite{G1} V 5).} 

On general grounds, the association of a localic groupoid to the set
of points of a topos 
is evident by means of an enrichment over localic spaces of the
categories of set-valued functors. Localic spaces are formal duals of
locales, and it is not evident how this enrichment
can be made in a way that furnish a  manageable theory for the sometimes 
unavoidable  work in the category of locales.  
Generalizing the construction in \cite{D} of the localic group of
automorphisms of a set-valued functor
we develop this enrichment in section \ref {enrichmentsection}, and  
in section \ref{lpoints} we construct the localic groupoid of
points. The objects are the points of the topos. The  hom-sets are (in
general pointless) localic spaces. This construction is
adequate for the representation theorems only in presence of enough
points in the topos. 

We develop in
detail the pointed theory in section \ref{Grothendieck's galois
theory}, where we bring into consideration the localic groupoid of all
the points. We establish the fundamental theorem in the form of a
characterization of
 Galois topoi with (at least one, and thus enough)
points as the classifying topoi of \emph{connected groupoids with
discrete space of objects and prodiscrete localic spaces of
hom-sets}. We also introduce the concept of \emph{proessential
point}, show how to construct Galois closures with this, and 
prove a new characterization of pointed Galois topoi. 

Grothendieck and Verdier always assume the existence
of enough points arguing that in all the meaningful examples the points are
there. Their thoughts on pointless topoi are revealed in
\cite{G2} Ex IV 6.4.2 where they write:  ``on peut
cependant, ``en faisant expres'' construire des topos qui n'ont pas
suffisamment de points''. However, with present hindsight, and as it was first
and long ago stressed by Joyal, we can argue that unpointed theories
are justified. 

 The theories in \cite{G2} and \cite{M2} are
localic only at the level of the fundamental groupoid
arrows. Fundamental groupoids of 
Galois topoi \emph{loose their objects by the same reason that they loose
their arrows} (namely, some co-filtered inverse limits of sets become
empty, see section \ref{unpointedsection}). It seems
natural then to develop a theory which is localic also at the level of
objects. 

M. Bunge in \cite{B} (see also
\cite{BM}) develop an unpointed theory for Galois topoi
following the inverse limit techniques implicit in \cite{G1} and
\cite{G2} and made explicit in \cite{M2}. Around the same time,
J. Kennison \cite{K} also developed an unpointed theory with a
different approach. They both prove the  fundamental theorem under 
the form of a Galois topoi characterization as the 
classifying topoi of \emph{prodiscrete localic groupoids}. 

Joyal-Tierney galois theory (see below) is behind M. Bunge
development of the
unpointed theory of Galois topoi. However, this theory follows by
inverse limit techniques directly from the theory of classifying topoi of
discrete groups (or groupoids), which is a very simple and elementary
case of Joyal-Tierney theorems. We show  in section
\ref{unpointedsection}
how the pointed theory of Galois topoi \emph{can as
well be developed} in an unpointed way along the same lines of 
\cite{M2}, \cite{B}, and \cite{BM}. We show that the localic groupoid
in the fundamental
theorem, even in the unpointed case, can be considered to be the
groupoid of (may be phantom) points of the topos.

The unpointed theory also applies in the presence of points, but it
yields an slightly different groupoid than the pointed theory. We
compare these groupoids in section \ref{comparisonsection}.

In their seminal paper  on galois
theory, \cite{JT} (after Grothendieck's \cite{G1}), Joyal and Tierney
bring new light into the subject. \emph{Galois theory is conceived by
  interpreting the fundamental theorem  
as an statement that says that a given geometric morphism of topoi
 is of effective descent} (namely, the point involved in the
classical and Grothendieck galois theories). They prove that \emph{any
  open surjection is of effective descent.} It follows an unpointed
theory of representation for a completely arbitrary topos in terms of
localic groupoids, which culminates
with their fundamental theorem \mbox{Ch.VIII 3. Theorem 2}, which
states that any topos is 
the classifying topos of a localic groupoid.  This theorem needs the
construction of a localic cover and sophisticated change of base
techniques, and we think it  
describes different phenomena from the one that concerns Galois topoi,
\emph{either pointed or unpointed}. The reader interested in
Joyal-Tierney theory of classifying topoi of localic groupoids should
also consult \cite{M3}.

The representation theorem of pointed connected atomic topoi,
\cite{JT}, Ch.VIII 3. Theorem 1, is  however closely related to the
representation theorem of pointed Galois topoi and classical galois
theory. It follows because any point of a connected atomic topos is an
open surjection, thus a geometric morphism of effective descent.   

In \cite{D} we
developed what we call \emph{localic galois theory} and prove therein
this result in a closer manner to classical galois theory,
independently of descent techniques (and of Grothendieck's inverse
limit techniques as well).  
This theorem shows that pointed connected atomic topoi classify
\emph{connected localic groupoids with discrete space of objects}. 
The groupoid in the theorem, as it is the case for Galois topoi, 
is the localic groupoid of points. We recall all this in section
\ref{localicgalois}.

Of course, Theorem 2 (loc. cit.) applies to an  
arbitrary, (even connected but maybe pointless) atomic topos, but
it is a different theorem. The localic groupoid  is not
canonically associated and cannot be considered to be 
(as far as we can imagine) a groupoid of (may be phantom) points of the
topos. Furthermore, when applied to atomic topoi with enough points
(one for each connected component suffices) it does not yield the
localic groupoid of points. 
It would be interesting to have a theorem  which, in
presence of a point, yields \mbox{Theorem 1.} We still
not know how to define the groupoid of phantom points for a general
atomic topos (as we do for a general Galois topos). 
An \emph{unpointed localic galois theory} (compressing the unpointed
prodiscrete galois theory as well as the pointed localic galois
theory) is yet to be developed.  

    

\tableofcontents

\enlargethispage*{1000pt}

{\sc acknowledgments.} 

I am grateful to M. Bunge for many fruitful and motivating discussions
on Galois topoi and its points or lack of them, and to A. Joyal
whom  in several occasions made me the gift of his deep insight on
Atomic and Galois topoi.

\pagebreak

\section{Background, terminology and notation} \label{background}

In this section we recall some topos and locale theory that we shall 
explicitly need, and in this way fix notation and terminology. We also
include some inedit proofs when it seems 
necessary. Our terminology concerning spaces and locales follows
Joyal-Tierney \cite{JT}, except that we define \emph{localic space} to
be the formal dual of a \emph{locale}, although we omit very often the
qualification ``localic'' and just write ``space''. Instead of saying
\emph{spatial group} we say \emph{localic group}, and the same for
groupoids. We do not distinguish notationally a localic space from its
corresponding locale. 

We denote $\cc{S}$ the topos of sets, and all topoi are
supposed to be Grothendieck topoi over $\cc{S}$. We nonetheless think
that all results in this paper hold as well for $\cc{S}$ an
arbitrary base Grothendieck topos, albeit, a few of them suitably
reformulated to avoid the use of choice.

\vspace{2ex} 

\subsection{Filtered inverse limits of topoi} \indent \vspace{1ex}
\label{inverselimitsection}

We recall here the fundamental result on filtered inverse limits of
topoi, which consists on the construction of the site for such kind of
limit.
Inverse limits of topoi have been extensively considered in SGA4, VI,
 where a fully detailed 2-categorical treatment is
developed. Consider a filtered system
of sites and morphisms of sites (continuous flat functors) and the
induced system of topoi as
shown in the following diagram (where the vertical arrows $\epsilon$
are the associate sheaf functor):
$$
\diagrama
 {
  \cc{C_{\alpha}} \ar[r]^{T_{\alpha \beta}} \ar[d]^{\epsilon}  
       &  \cc{C_{\beta}} \ar[d]^{\epsilon}  
       &  \cdots  \ar[r]  &  \cc{C}  \ar[d]^{\epsilon} 
       & & \cc{C_{\alpha}} \ar[r]^{T_{\alpha}} \ar[d]^{\epsilon} 
       &  \cc{C}  \ar[d]^{\epsilon} \\
  \cc{C_{\alpha}^{\sim}} \ar[r]^{t_{\alpha \beta}^{*}}  
       &  \cc{C_{\beta}^{\sim}}  
       &  \cdots  \ar[r]  &  \cc{C^{\sim}}
       & & \cc{C_{\alpha}^{\sim}} \ar[r]^{t_{\alpha}^{*}} 
       & \cc{C^{\sim}}
  }
$$
The diagram $\cc{C_{\alpha}} \mr{T_{\alpha}} \cc{C}$  is the filter colimit
of the categories $\cc{C_{\alpha}}$, and the category $\cc{C}$ is
furnished with the coarsest topology that makes the
inclusions $T_{\alpha}$ continuous. The resulting site is called 
\emph{the inverse limit site}. It is shown in \cite{G2} that the
inclusions are flat (here is where the filterness condition plays a
key role), and thus they are morphisms of sites.  With this at hand,
the next theorem follows immediately from SGA4, Ex. IV, 4.9.4 
.  

\begin{theorem} \label{inverselimittopoi} \emph{(SGA4, Ex VI,
    8.2.11)}. In the
situation described above, the following formula 
\mbox{$Lim_{\alpha} (\cc{C_{\alpha}})^{\sim} = (Colim_{\alpha} \, 
\cc{C_{\alpha}})^{\sim}$} holds. That is, in the diagram above, the bottom
row consists of the inverse image functors of a filtered inverse limit
of topoi.

\qed \end{theorem}

The interested reader will profit also consulting \cite{M1}, where
many ubiquitous and important preservation properties of filtered
inverse limits are stated and proved. There, a construction of the
inverse limit site (theorem 3.1 loc. cit.) is developed in the style
of the classical construction of the \emph{p-adic numbers}. It is
straightforward to  check that this construction (made for inverse
limit sequences) can be as easily done for general filtered systems, a
fact that has its own independent interest. Then, all results
in \cite{M1} can be derived directly for general filtered inverse
limits in the same way that for sequences.

\vspace{2ex}

\subsection{Basic facts on posets and locales} 
\indent \vspace{1ex} \label{localesection}

We think of locale theory as a 
reflection of topos theory (with the poset $2 = \{0, 1\}$ playing the 
role of the category $\cc{S}$ of sets), as well as that of a 
theory of generalized topological spaces.

\vspace{1ex}
 
We consider a \emph{poset} as a category, and in this vein a 
\emph{partial order} is a reflexive and transitive 
relation, not necessarily antisymmetric. We 
shall refer to the elements of a poset as \emph{objects}.

\vspace{1ex}

Given any poset $D$, the free inf-lattice on $D$, which we denote
$\cc{D}(D)$, is furnished with a poset-morphism $D \mr{\eta}\cc{D}(D)$
which  is generic in the sense that giving any inf-lattice $H$,
composing 
with $\eta$ defines an equivalence of posets
\mbox{$\mathcal{L}ex(\cc{D}(D),\, H) \mr{\simeq} 
\mathcal{P}os(D, H)$}, where  $\mathcal{L}ex(\cc{D}(D),\, H)$ and 
$\mathcal{P}os(D, H)$ indicate
inf-preserving morphisms and poset-morphisms respectively.

We recall now a construction of $\cc{D}(D)$. 

\begin{proposition} \label{freeinflattice}
The objects of $\cc{D}(D)$ are in one to one correspondence with the
finite subsets of $D$. Given a subset 
$A = \{a_{1}, \ldots ,\,a_{n}\} \subset D$, we denote 
\mbox{$[A] \;=\; [\langle a_{1}\rangle , \ldots ,\,\langle a_{n}\rangle ]$} the corresponding
object in $\cc{D}(D)$. The morphism $\eta$ is defined by $\eta(a) =
[\langle a\rangle ]$. Given any other \mbox{$[B] = \,[\langle b_{1}\rangle , \ldots 
,\,\langle b_{k}\rangle ]$}, the order relation is given by:
$$
  \frac
       {
        \;\;\;
        [\langle a_{1}\rangle , \ldots ,\,\langle a_{n}\rangle ]
  \;\;\leq \;\;[\langle b_{1}\rangle , \ldots ,\,\langle b_{k}\rangle ]
        \;\;\;
       }
       {
        \;\exists\:\sigma\,:\,  \{1, \ldots, k\} 
                          \mr{} \{1, \ldots, n\},\; a_{\sigma i} \leq b_{i}\,
       }
$$
\qed  \end{proposition}

A \emph{locale} is a complete lattice in which finite infima 
distribute over arbitrary suprema. A morphism of locales 
$E \mr{f^{*}} H$ 
is defined as a function $f^{*}$ preserving finite infima and 
arbitrary suprema (notice that we put automatically an upper star 
to indicate that these arrows are to be considered as 
inverse images of geometric maps).

\vspace{1ex}

 \emph{Inf-lattices} $D$ are sites of definition for locales (rather 
than bases 
 of opens). \mbox{2-valued} \emph{presheaves} $D^{op} \rightarrow 2$ 
 form a locale, $D^{\wedge} = 2^{D^{op}}$. 
Given a 
  Grothendieck 
 (pre) topology on $D$, 2-valued \emph{sheaves} also 
form a 
 locale, denoted $D^{\sim}$. The associated sheaf defines a morphism 
of locales 
  \mbox{$D^{\wedge} \rightarrow  D^{\sim}$,} and this is a procedure 
in which 
  quotients of locales are obtained. A \emph{site} is, in this sense, a 
  \emph{presentation} of the locale of sheaves.

 The basic fundamental result of this construction is the following:

\begin{lemma} \label{diaconescu}
The associated sheaf
\mbox{$D \mr{\#} D^{\sim}$} is a morphism of sites (preserves infima
and sends covers into epimorphic families) into a locale 
which is generic, in the sense 
that giving any locale $H$, composing with $\#$ defines an 
equivalence of posets
\mbox{$\mathcal{M}orph(D^{\sim}, H) \mr{\simeq} \mathcal{M}orph(D, 
H)$}.
\qed \end{lemma}

This lemma is just \cite{G2} IV 4.9.4 in the poset context.

\vspace{1ex}

 A \emph{localic space} is the formal 
dual of a locale. Thus, $E \mr{f^{*}} H$ defines a map or morphism of 
localic spaces from $H$ to $E$, $H \mr{f} E$. Following \cite{JT}, 
all 
these maps are called \emph{continuous maps}.  

\vspace{1ex}

The \emph{open subspaces} of a localic space E correspond to the
objects of the locale E (\cite{JT} ch V, 2.). We shall identify 
(as an abuse of notation)
the object $u \in E$ with the subspace defined by the quotient
locale $E \mr{} U$, $w \mapsto w \wedge u$, 
$U = \{v | v \leq u\}$. We abuse $u = U$ and  indistinctly write 
$u \subset E$ or $u \in E$.

\vspace{1ex}
 
A \emph{surjection} between localic spaces is a map whose inverse
image reflects isomorphisms. It follows immediately from the
preservation of infima that $f^{*}$ is injective (up to
isomorphisms). Thus, surjections are epimorphisms in the category of
localic spaces. Furthermore, it also follows that $f^{*}$
is full, in the sense that the implication 
$(f^{*}u  \leq f^{*}v \; \Rightarrow \;  u \leq v)$ holds.

\vspace{1ex}

A \emph{localic monoid}, (resp. \emph{localic group}) is a monoid 
object (resp. group object) in the category of localic spaces. 
A \emph{morphism of monoids (or groups)} $H \mr{\varphi} G$ is a 
continuous 
map satisfying the usual identities. Actually, all this is given in 
practice by the inverse image maps between the corresponding locales 
satisfying the dual equations.
  
\vspace{1ex}

The \emph{locale of relations} $lRel(X, \, Y)$ between two sets $X$, $Y$  is
the  free locale on 
$X \times Y$.  Recall that the free locale on a set $S$ is
constructed by taking presheaves on the free inf-lattice on $S$ (the
lattice of finite subsets  with the dual order, see
\ref{freeinflattice}). 
If $\{(x_{1},\,y_{1})\,\ldots , (x_{n},\,y_{n})\} \subset X \times Y$, 
we write 
\mbox{$[\langle x_{1}\,|\,y_{1}\rangle \,\ldots, \langle x_{n}\,|\,y_{n}\rangle ]$} for the 
corresponding object in the inf-lattice and in the locale. Remark that
this object is the finite infimun of the $(x_{i},\,y_{i})$ 
(see \cite{D} for details, there for the case $X = Y$). 

\vspace{1ex}

G. Wraith in an inspiring paper \cite{W} defines the locales of
functions and of bijections between two sets $X$ and $Y$ by
considering the appropriate  generators and relations. In our context
these relations become covers in the free inf-lattice on $X \times Y$.    

\vspace{1ex}

The \emph{locale of functions} $lFunc(X, \, Y)$ from $X$ to $Y$, is the locale
of sheaves for the topology that forces a relation to be a function 
(\cite{W}, \cite{D}). It is generated by the following covers (u: univalued)
and (e: everywere defined).:
$$
\begin{array}{l} 
 u)\;\; \emptyset \;\rightarrow\; [\langle \!z\,|\,x\!\rangle ,\;\langle \!z\,|\,y\!\rangle ] 
 \;\;\;\;\;\;\;\;\;\;\;
 (each \; z \in X \,, \; x \neq y \in Y)  
 \\[1ex]
 e)\;\; [\langle \!z\,|\,x\!\rangle ] \;\rightarrow \; 1, \; x \in Y 
 \;\;\;\;\;\;\;\;\;\;\;\;\;\;\;
 (each \; z \in X)
\end{array}
$$

\vspace{1ex}

The \emph{locale of bijections} $lBij(X, \, Y)$ is
determined if we add the covers which force a function to be bijective
(i: injective) and (s: surjective):
$$
\begin{array}{l} 
 i)\;\; \emptyset \;\rightarrow\; [\langle \!x\,|\,z\!\rangle ,\;\langle \!y\,|\,z\!\rangle ] 
 \;\;\;\;\;\;\;\;\;\;\;
 (each \;  x \neq y \in X \,, \;  z \in Y)  
 \\[1ex]
 s)\;\; [\langle \!x\,|\,z\!\rangle ] \;\rightarrow \; 1, \; x \in X 
 \;\;\;\;\;\;\;\;\;\;\;\;\;\;
 (each \; z \in Y)
\end{array}
$$

We denote $lAut(X) = lBij(X, \, X)$.

\vspace{1ex}

Given any sets $X$, $Y$, corresponding to the basic covers the following
equations hold in the locale $lFunc(X, \, Y)$:
$$
u)\;[\langle \!z\,|\,x\!\rangle ,\;\langle \!z\,|\,y\!\rangle ] \;=\; 0 \;\;\; 
(each \; x \neq y,\; z)\,, \;\;\;\;
e)\;\bigvee\nolimits_{x} \; [\langle \!z\,|\,x\!\rangle ] \;=\; 1 \;\;\;
(each \, z)
$$
Two additional equations hold in $lBij(X, \, Y)$:
$$
i)\;[\langle \!x\,|\,z\!\rangle ,\;\langle \!y\,|\,z\!\rangle ] \;=\; 0 \;\;\; 
(each \; x \neq y,\; z)\,, \;\;\;\;
s)\;\bigvee\nolimits_{x} \; [\langle \!x\,|\,z\!\rangle ] \;=\; 1 \;\;\;
(each \, z)
$$  

Notice that 
we abuse notation and omit to indicate the associated sheaf morphism. 

Given any set $X$, we consider now the group structure on the localic
space $lAut(X)$. More generally, 
it is tedious but straightforward to check the following:
\begin{proposition} \label{groupstructure}
For any sets  $X$, $Y$, $Z$, there are morphisms of locales:
$$
lFunc(X, \, Y)\;\; \mr{m^{*}}\;\; lFunc(X, \, Z) \otimes 
lFunc(Z, \,Y)\,,\ \;\;  lFunc(X, \, X)\;\; \mr{e^{*}}\;\;2 
$$
defined on the generators by the following formulae:
$$
m^{*}[\langle \!x \:|\:y\!\rangle ] \;=\;
 \bigvee\nolimits_{z} \;[\langle \!x \:|\:z\!\rangle ]\otimes[\langle \!z \:|\:y\!\rangle ]\,,\;\;
e^{*}[\langle \!x\,|\,y\!\rangle ] \;=\; 1 \;\;\Leftrightarrow\;\; x \;=\; y 
$$
These data satisfy the equations of an enrichment of the category of
sets $\cc{S}$ over the category
of localic spaces, which we shall denote $l\cc{S}$. In particular,
for each set $X$, the localic space 
\mbox{lFunc(X, X)} is a localic monoid.

The above formulae together with 
\mbox{$\iota^{*}[\langle \!x \:|\:y\!\rangle ] \;=\;[\langle \!y
    \:|\:x\!\rangle ]$} also define  morphisms of locales:
$$
lBij(X, \, Y) \;\;\mr{m^{*}} \;\;lBij(X, \, Z) \otimes lBij(Z, \, Y)
\,,\;\;\; lBij(X, \, X)\;\;\mr{e^{*}}\;\; 2
$$
$$
lBij(X, \, Y)\;\; \mr{\iota^{*}} \;\;lBij(Y, \, X)
$$
which determine an structure of localic groupoid on the (discrete) set
of all sets. In particular, for each set $X$, the localic space 
$lAut(X)$ is a localic group.   \qed   
\end{proposition}

\vspace{2ex}

\subsection{The classifying topos of a localic groupoid} \indent \vspace{1ex}

Following \cite{W} we now define group actions in terms of the
(sub base) generators of the localic group $lAut(X)$:

\begin{definition} \label{action}
Given a localic group $G$ and a set $X$, an \emph{action} of $G$ on 
$X$ is a continuous morphism of localic groups $G \mr{\mu} lAut(X)$. 
It 
is completely determined by the value of its inverse image on the 
generators, $X \times X  \mr{\mu^{*}} G$. By definition, the 
following 
equations hold
$$
m^{*}\mu^{*} \;=\; (\mu^{*} \otimes \mu^{*})\,m^{*}\,,\;\;\;\;\;\;
\mu^{*}\iota^{*} \;=\; \iota^{*}\mu^{*}\,,\;\;\;\;\;\; 
e^{*}\mu^{*} \;=\; e^{*}
$$
We say that the action is 
\emph{transitive} 
when for all \mbox{$x \in X,\: y \in X$, $\mu^{*}[\langle x\,|\,y\rangle ] \neq 
0$}.  
\end{definition}
Given a localic group $G$, a \emph{$G$-set} is a set furnished with 
an action of $G$. 
\begin{definition} \label{Gsetarrows}
Given two $G$-sets $X$, $Y$,
 $X \times X  \mr{\mu^{*}} G$, \mbox{$Y \times Y  \mr{\mu^{*}} G\,$,} 
a morphism 
 of $G$-sets is a function $X \mr{f} Y$ such as 
\mbox{$\mu^{*}[\langle \!x \:|\:y\!\rangle ] \;\leq\; \mu^{*}[\langle \!f(x) 
\:|\:f(y)\!\rangle ]$.} 
This defines a category $\mathcal{B}G$ furnished with an  
underlying set functor $\mathcal{B}G \mr{} \cc{S}$ 
into the category of sets.
\end{definition}
\begin{definition} \label{lFix(x)}
Given a localic group $G$ acting on a set $X$, and an element $x\in X$, 
the open subgroup of $G$, informally described
 as $\{g\in G \;|\; gx = x \}$, is defined to be the object 
 \mbox{$lFix(x) = \mu^{*}[\langle x\,|\,x\rangle ]$} in the locale $G$.
\end{definition}

 Given a morphism between two localic groups 
\mbox{$G \mr{t} H$}, and an action of $H$ in a set $X$, 
$X \times X \mr{\mu^{*}} H$, the composite
\mbox{$X \times X \mr{\mu^{*}} H \mr{t^{*}} G$} defines an action of
$G$ on $X$. This defines a functor, that we denote $\cc{B}(t)^{*}$,
$\mathcal{B}H \mr{} \mathcal{B}G$ (clearly commuting with the
underlying sets), and all these assignments are
functorial in the appropriate sense. 

The transitive $G$-sets are the connected objects of $\mathcal{B}G$. 
We shall denote 
$t\mathcal{B}G$ the full subcategory of non empty transitive $G$-sets.

\begin{proposition} \label{gsetatomic}
The category $t\mathcal{B}G$ is an small category, which together with
the underlying set functor satisfies 1) i) ii) iii) iv) in proposition
\ref{atomicsite}. The topos of sheaves for the canonical topology is 
$\mathcal{B}G$, which is then a pointed connected atomic topos. 
The canonical point, that we denote $u$, has the inverse
image given by the underlying set functor. 
\end{proposition}
\begin{proof}
This is proposition 8.2 in \cite{D}.
\end{proof}

\begin{proposition} \label{surjections}
Given a morphism of localic
groups \mbox{$G \mr{t} H$,} the functor 
$\mathcal{B}H \mr{\cc{B}(t)^{*}} \mathcal{B}G$ is the inverse image of a
morphism of pointed topoi. If the morphism is a surjection, then,
given any transitive $H$-set $X$, 
$\cc{B}(t)^{*}(X)$ is a transitive $G$-set, and the functor 
$t\mathcal{B}H \mr{\cc{B}(t)^{*}} t\mathcal{B}G$ is full and faithful.      
\end{proposition}
\begin{proof}
The first assertion is
   straightforward and stated in, for example, \cite{M2}. The second
   assertion is immediate: $0 \langle  \mu^{*}[\langle \!x \:|\:y\!\rangle ]$, then 
$t^{*}\mu^{*}[\langle \!x \:|\:y\!\rangle ]$
   cannot be equal to $0$ since $t^{*}$ reflects isomorphisms by
   definition. Finally, let $X$, $Y$ be any two $H$-sets and 
   $X \mr{f} Y$ a morphism for the $G$-actions, that is,
   \mbox{$t^{*} \mu^{*}[\langle \!x \:|\:y\!\rangle ] \leq t^{*} \mu^{*}[\langle \!fx
   \:|\:fy\!\rangle ]$.} 
Since inverse
   images of surjections between locales are full, it follows that 
   \mbox{$\mu^{*}[\langle \!x \:|\:y\!\rangle ] \leq \mu^{*}[\langle \!fx \:|\:fy\!\rangle ]$}.  
\end{proof}

\vspace{2ex}
 
Given any localic groupoid $\cc{G}$,  the category of \emph{discrete
$\cc{G}$ spaces} 
is defined in a standard way in \cite{JT} VIII, 3, and proved
therein to be a topos, denoted $\cc{B}\cc{G}$ (see also \cite{M3} 5.2).

Consider the enrichment $l\cc{S}$ of the category of sets over the
category of localic spaces, \ref{groupstructure}. It is straightforward
to check the following:

\begin{proposition} \label{betadiscreteobjects}
Given any localic groupoid with discrete set of objects, the category
$\cc{B}\cc{G}$  can be
defined as the (ordinary) category of enriched functors 
$\cc{G} \mr{} l\cc{S}$ and natural transformations. 
$\cc{B}\cc{G} =  l\cc{S}^{\cc{G}}$. In turn, it is
straightforward to define these data in the style of definitions
\ref{action} and \ref{Gsetarrows}. For each object of the
groupoid there is a corresponding evaluation functor 
$\cc{B}\cc{G} \mr{} \cc{S}$, 
and these functors (collectively) reflect isomorphisms.
\qed \end{proposition}

A localic groupoid with discrete set of objects is said to be
\emph{connected} if for
each pair of objects $p, \; q$, the localic space $\cc{G}[p,\, q]$ is
non empty (equivalently, in the notation of \cite{M3}, if the morphism 
$G_1 \mr{(d_0, \, d_1)} G_0 \times G_0$ is a surjection). A connected
localic groupoid may not be connected as an ordinary groupoid
since the localic spaces $\cc{G}[p,\, q]$ can be pointless 
(see example \ref{example1}).

It is possible to check with the methods of
\cite{D} that the tops $\cc{B}\cc{G}$ is atomic, 
and that if the groupoid is connected,
it is equivalent to $\cc{B}G_p$, where 
$G_p = \cc{G}[p, \,p]$ is any one of its vertex localic
groups (notice that the first  assertion follows from the second
(and \ref{gsetatomic}), since a sum of atomic topoi is atomic). 
We omit to do all this in print, and invoke \cite{M3} as a proof. 

\begin{proposition} \label{gsetatomic1}
Given any localic groupoid $\cc{G}$ with discrete set of objects, the
topos $\cc{B}\cc{G}$ is an atomic topos with enough points (with
inverse images given by the evaluation functors). If the groupoid is
connected, then it is a connected topos, equivalent to the
classifying topos of any one of its vertex localic groups.
\end{proposition}
\begin{proof}
It is stated in \cite{M3}, 4.5 c) that it is an atomic topos. Clearly, it
has enough points. The second statement follows immediately from
\cite{M3} 5.15 (v) considering the inclusion morphism. 
\end{proof}

\vspace{2ex}
 
\subsection{Classifying topos and filtered inverse limits} \indent
\vspace{1ex}  
\label{locinvsection}

In this subsection we study the behavior of the classifying topos
regarding filtered inverse limits of localic groups, see \cite{M2}, and
groupoids. The reader should be aware that for the
applications to the representation of Galois topoi, only the
particular case of localic limits of discrete groups (or groupoids) is
necessary.

 Consider a localic group $G$ an open subgroup $u \subset  G$. 
In \cite{M2} 1.2, the quotient
 localic space \mbox{$G \mr{\rho} G/u$} is defined as usual in group theory,
 working formally in in the category of localic spaces considered as
the formal dual of the category of locales. Then, it is proved
 that it is a discrete localic space, and that the localic group
 $G$ has an ``obvious'' transitive action in the set $Z = [G, \: 2]$
 of its points (\cite{M2} 2.3). Furthermore, $u = lFix(z_{0})$ for the point 
$z_{0} \in Z$ defined by the composite $e^{*} \circ \rho^{*}$ (where 
$G \mr{e^{*}} 2$ is the (co) unit of $G$). It
 follows from \cite{D} proposition 7.9 that any transitive action is
 of this form. 

It is also stated in \cite{M2} 2.4 that given any (co) filtered inverse
limit of localic groups $G_{\alpha} \ml{t_{\alpha}} G$, the subgroups
of $G$ of the form $t^{*}_{\alpha}(w)$ for some open subgroup 
$w \subset G_{\alpha}$ form a cofinal system of open subgroups of $G$, in
the sense that given any open subgroup $u \subset G$, there exists
$\alpha$ and an open subgroup $w \subset G_{\alpha}$ such that 
$t^{*}_{\alpha}(w) \leq u$ (this is due to the fact that the objects
of $G$ of the form $t^{*}_{\alpha}(w)$ for some $\alpha$ and 
$w \in G_{\alpha}$ generate $G$, as it can be seen, for example, by 
theorem \ref{inverselimittopoi} in the context of posets and locales).      

We prove now a generalization of a classical result in the theory of
profinite topological groups.

\begin{proposition} \label{factoring}
Consider a (co) filtered inverse limit diagram $G_{\alpha} \ml{t_{\alpha}} G$
of localic groups with
surjective transition morphisms $G_{\beta} \ml{t_{\alpha \beta}} G_{\beta}$.
Then, given any transitive $G$-set X,  \mbox{$X \times X \mr{\mu^{*}} G$}, 
the action $\mu$ 
factors through some $G_{\alpha}$-action in the following sense:

There exists $\alpha$, a $G_{\alpha}$-set $Z$, 
$Z \times Z \mr{\mu^{*}} G_{\alpha}$,  an epimorphism of $G$-sets 
$\cc{B}(t_{\alpha})^{*}(Z) \mr{f} X$ (given by a  surjective function 
$Z \mr{f} X$), and a factorization as follows:
$$
\diagrama
         {
             Z \times Z  \ar[d]^{\mu^{*}}  \ar[r]^{f \times f}  
          &  X \times X \ar[d]^{\mu^{*}} 
          \\
          G_{\alpha} \ar[r]^{t^{*}}  &  G \ar@{}[lu]|{\le}
         }
$$        
\end{proposition}
\begin{proof}
Choose any $x_{0} \in X$, and consider the open subgroup 
\mbox{$lFix(x_{0}) = \mu^{*}[\langle \!x_{0} \:|\:x_{0}\!\rangle ] \in G$.} By the
remarks preceding this proposition, 
there exists $\alpha$, a $G_{\alpha}$-set $Z$, 
\mbox{$Z \times Z \mr{\mu^{*}} G_{\alpha}$,} and an element $z_{0} \in Z$ such
that \mbox{$t_{\alpha}^{*}lFix(z_{0}) = 
t_{\alpha}^{*} \mu^{*}[\langle \!z_{0} \:|\:z_{0}\!\rangle ] \leq lFix(x_{0})$.} 
But $t_{\alpha}^{*}\mu^{*}[\langle \!z_{0} \:|\:z_{0}\!\rangle ]$ is the
subgroup $lFix(z_{0})$ for the action
$\cc{B}(t_{\alpha})^{*}(Z)$. Since the projection $t_{\alpha}$ is
surjective (see \cite{JT}), this action is transitive (cf \ref{surjections}).
 The proof finishes then by \cite{D}, 7.9.   
\end{proof}

We can improve now a little over \cite{M2}, where the following
theorem is proved in the case of open surjections.

\begin{theorem} \label{inverselimitgroups}
Given a (co)filtered diagram of localic groups and 
surjective localic group morphisms, and its inverse limit:
$$ G_{\alpha} \ml{t_\alpha \beta} G_{\beta} \;\;\;\;\; \dots  \ml{} G$$
the induced diagram of topoi and topoi morphisms:
$$ \cc{B}(G_{\alpha}) \ml{\cc{B}(t_{\alpha \beta})}   
   \cc{B}(G_{\beta}) \;\;\;\;\; \dots  \ml{} \cc{B}G$$
is also an inverse limit diagram.
\end{theorem}
\begin{proof}
We have the situation described in the following diagram;
$$
\diagrama
      {
       t\cc{B}G_{\alpha}\; \ar@{^{(}->}[r]^{\;\;\cc{B}(t_{\alpha
                                              \beta})^{*}} \ar[d]^{\epsilon}  
       & \;t\cc{B}G_{\beta}\;  \ar[d]^{\epsilon}  
       &  \cdots  \ar@{^{(}->}[r]  
       &  \;\cc{C}\;\;  \ar[d]^{\epsilon} \ar@{^{(}->}[r] 
       & \;t\cc{B}G  \ar[d]^{\epsilon} 
       \\
       \cc{B}G_{\alpha}\; \ar[r]^{\cc{B}(t_{\alpha \beta})^{*}}  
       & \;\cc{B}G_{\beta}\;    
       &  \cdots  \ar[r]  
       &  \;\cc{C}^{\sim}\;  \ar[r] 
       & \;\cc{B}G 
      }
$$
Here $\cc{C}$ is the inverse limit site (as a category $\cc{C}$ is
the filtered colimit), and in the top row all functors are full and
faithful by proposition \ref{surjections}. By \ref{inverselimittopoi}
above $\cc{C}^{\sim}$ is
the inverse limit of the topoi $\cc{B}G_{\alpha}$. Notice that the
topology in $\cc{C}$ is induced by the canonical topology of 
$t\cc{B}G$. Then, by \ref{factoring}, it follows from the comparison
lemma (\cite{G2}, Expos\'e III, 4) that the arrow 
$\cc{C}^{\sim} \mr{} \cc{B}G$ is an equivalence.
\end{proof}

We comment that the corresponding theorem for filtered inverse limits
of \emph{discrete groupoids} has been stated and proved with do care by Kennison
 \cite{K} 4.18. In the result's statement it is necessary to
assume that transition morphisms are \emph{composably onto}, (see
\cite{K}).  This takes care of the necessary surjectivity of the
system at the level of arrows. In view of \ref{betadiscreteobjects}, the second statement
in \ref{gsetatomic1}, and
\ref{inverselimitgroups} above, a similar theorem for filtered inverse
limits of localic groupoids with discrete sets of objects seems plausible.  
Abusing rigor, one could say that a corresponding result in the
case of arbitrary localic groupoids also holds, but we do not know of
any clear proof in print, and it remains an \emph{open problem}
to us.

\vspace{2ex} 

\section{Enrichment of set valued functor categories over
localic spaces}  

In this section we do a brief review of the salient features of the
construction and properties given in \cite{D} of the locale of
automorphism of a set-valued functor. We develop the more
general case of natural transformations between two functors. We
establish the whole 2-categorical specifications in some cases, and
prove some new results. We introduce the localic groupoid of points of
a topos, and study its behavior regarding filtered inverse limits. It is
pertinent to remark that only the particular case of a system where
the points are representable (thus the groupoids discrete, but not
their limits) is necessary for the applications to the representation
of Galois topoi.

\vspace{2ex}

\subsection{The localic functor category} \indent \vspace{1ex}
\label{enrichmentsection}

Given any category $\mathcal{C}$ and any set-valued functor
$F:\mathcal{C} \mr{} \cc{S}$, recall that the diagram of $F$, 
which we denote $\Gamma_{F}$, is the category whose objects are
the elements of the disjoint union of the sets $FX, \; X \in 
\mathcal{C}$. That is to say, pairs $(x,X)$ where $x \in FX$. The arrows 
$(x,X) \mr{f} (y,Y)$ are maps $X \mr{f} Y$ such that $F(f)(x)=y$.
Given any $(z, Z)$ in $\Gamma_{F}$, there is a natural transformation 
$[Z, \:-] \mr{z^*} F$ defined as follows: given $h \in [Z, \:X], z^{*}(h) =
F(h)(z)$, and the resulting diagram is a colimit cone 
(indexed by  $\Gamma_{F}$). 

Associated with $\Gamma_{F}$, we define a poset, which we denote $D_{F}$,
identifying all arrows in each hom-set of category $\Gamma_{F}$.

\vspace{1ex}
 
Given any category $\mathcal{C}$ and any pair of set-valued functors
\mbox{$F:\mathcal{C} \mr{} \cc{S}$,} 
\mbox{$G:\mathcal{C} \mr{} \cc{S}$,}
a \emph{natural relation} between $F$ and $G$ is a 
relation $R \subset F \times G$ in the functor category. That is, it 
is a family of relations $RX$ on $FX \times GX$, $X \in \mathcal{C}$, 
such that given 
any arrow $X \mr{f} Y$ in $\mathcal{C}$, 
$(Ff \times Gf)RX \;\subset\; RY$. 
In other terms, it is a family of functions $FX \times GX
\mr{\phi_{X}} 2$ such that 
$\phi_{X} \:\leq \: \phi_{Y} \circ (Ff \times Gf)$.
It is clear that if a natural relation is functional, then it is a 
natural transformation.

In \cite{D} the locale of natural 
relations  from $F$ to $G$ is constructed and characterized as follows:

Consider the composite of the diagonal functor 
 $\mathcal{C} \rightarrow \mathcal{C} \times \mathcal{C}$
  with $F \times G$, which we denote $F \Delta G$, 
$(F \Delta G)(X) \,=\, FX \times GX$.
Consider the poset $D_{F \Delta G}$ whose objects are the disjoint 
union of the sets \mbox{$FX \times GX,\; X \in \mathcal{C}$.} 
The order relation is given by the following rule:
$$
  \frac{(X,\: (x_{0},\,x_{1})) \;\;\leq\;\; (Y,\: (y_{0},\,y_{1}))}
  {\; \exists \; X \mr{f} Y \; F(f)(x_{0}) = y_{0}\,, \; G(f)(x_{1}) = y_{1}}
$$
 
Consider then the free inf-lattice $\mathcal{D}(D_{F \Delta G})$ 
on this poset (see \ref{freeinflattice}). The locale of presheaves on
this lattice is the locale of natural relations from $F$ to $G$. By
introducing in $\mathcal{D}(D_{\Delta FG})$ the appropriate covers,
we construct the (quotient) locales of natural transformations and
natural bijections.   

Given an object $X$, a pair $(x_{0},\,x_{1}) \in FX \times GX$ and a
 finite subset $A \subset D_{F \Delta G}$, we denote  
$$[(X,\,\langle \!x_{0}\,|\,x_{1}\!\rangle ),\; A]
\;=\; [(X,\,\langle \!x_{0}\,|\,x_{1}\!\rangle )]\,\wedge\, [A]$$ 
the corresponding object in $\mathcal{D}(D_{F \Delta G})$.


For each $X \in \mathcal{C}$ there is a function
$FX \times GX \mr{\lambda_{X}} \mathcal{D}(D_{F \Delta G})$
defined by $\lambda_{X}(x_{0},\,x_{1}) = 
[(X,\,\langle \!x_{0}\,|\,x_{1}\!\rangle )]$.

\begin{proposition} \label{genericnr} $\,$ \\
\indent 1) 

1.1)The locale 
\mbox{$lRel(F, \,G) \,=\, \mathcal{D}(D_{\Delta FG})^{\wedge}$} 
of natural relations from $F$ to $G$ is the  
locale of presheaves on  $\mathcal{D}(D_{\Delta FG})$
(we consider this inf-lattice as a site with the empty topology).

1.2) The locale 
\mbox{$lFunc(F, \,G) \,=\, \mathcal{D}(D_{\Delta FG})^{\sim}$}
of natural transformations from $F$ to $G$ is the  
locale of sheaves for the topology on  $\mathcal{D}(D_{\Delta FG})$ 
which forces a natural relation to be functional.   
This topology is generated by the following basic covers: (u: univalued)
and (e: everywere defined).
$$
\begin{array}{l} 
 u)\;\; \emptyset \;\rightarrow\; [(X,\;\langle \!z\,|\,x\!\rangle ),\;(X,\;\langle \!z\,|\,y\!\rangle )] 
 \;\;\;
 (each\; X,\; and \; each \; z \in FX \,, \; x \neq y \in GX)  
 \\[1ex]
 e)\;\; [(X,\;\langle \!z\,|\,x\!\rangle )] \;\rightarrow \; 1, \; x \in GX 
 \;\;\;\;\;\;\;\;\;\;\;
 (each \;X \;and \;each \; z \in FX)
\end{array}
$$
\indent 1.3) The locale  
\mbox{$lBij(F, \,G) \,=\, \mathcal{D}(D_{\Delta FG})^{\sim}$} of
natural bijections is constructed if we add the following covers:
(i: injective) and (s: surjective)
$$
\begin{array}{l} 
 i)\;\; \emptyset \;\rightarrow\; [(X,\;\langle \!x\,|\,z\!\rangle ),\;(X,\;\langle \!y\,|\,z\!\rangle )] 
 \;\;\;
 (each\; X,\; and \; each \;  x \neq y \in FX \,, \;  z \in GX)  
 \\[1ex]
 s)\;\; [(X,\;\langle \!x\,|\,z\!\rangle )] \;\rightarrow \; 1, \; x \in FX 
 \;\;\;\;\;\;\;\;\;\;\;
 (each \;X \;and \;each \; z \in GX)
\end{array}
$$

The
points of these locales are exactly natural relations, natural
transformations, and  natural bijections respectively.

\indent 2)

2.1) The inf-lattice $H = \mathcal{D}(D_{\Delta F})$, together with
     the functions $\phi_X = \lambda_X$ satisfy the
     following condition: 
 
 For each $X \in \mathcal{C}$, there is a function 
$FX \times GX \mr{\phi_{X}} \mathcal{D}(D_{F \Delta G})$, 
such that for each $X \mr{f} Y$,  
$\phi_{X} \leq \phi_{Y} \circ (F(f) \times G(f))$. 


\vspace{1ex}

2.2) The site defined in 1.2 satisfy in addition,

The following families are coverings: 
$$
\begin{array}{l} 
 \;\;u)\;\;  \emptyset \;\rightarrow\; \phi_{X}(z, x) \wedge \phi_{X}(z, y) 
 \;\;\;\;\;\;\;
 (each\; X,\; and \; each \; z \in FX \,, \; x \neq y \in GX)  
 \\[1ex]
 \;\;e)\;\; \phi_{X}(z, x) \;\rightarrow \; 1, \; x \in GX 
 \;\;\;\;\;\;\;\;\;
 (each \;X \;and \;each \; z \in FX)
\end{array}
$$
\indent 2.3) The site defined in 1.3  satisfy in addition,

The following families are also coverings:
$$
\begin{array}{l} 
 \;\;i)\;\; \emptyset \;\rightarrow\; \phi_{X}(x, z) \wedge \phi_{X}(y, z) 
 \;\;\;\;\;\;\;
 (each\; X,\; and \; each \;  x \neq y \in FX \,, \;  z \in GX)  
 \\[1ex]
 \;\;s)\;\; \phi_{X}(x, z) \;\rightarrow \; 1, \; x \in FX 
 \;\;\;\;\;\;\;\;\;
 (each \;X \;and \;each \; z \in GX)
\end{array}
$$            
\indent 3) These sites have, and therefore are characterized by, 
the following universal property:           
 
For any other 
such data, $FX \times GX \mr{\phi_{X}} H$, there is a unique morphism
of sites $\phi$ (as indicated in the diagram below): 
$$
\diagrama
        {
         FX \times GX \ar[rd]^{\lambda_{X}}  \ar@(r, ul)[rrd]^{\phi_{X}} 
                                                \ar[dd]_{F(f) \times G(f)}  
         \\
         \ar@{}[r]^{\leq\;\;\;\;\;\;\;\;\;\;\;\;\;\;}
         & 
         \;\;\cc{D}(D_{\Delta FG})\;\; \ar@{-->}[r]^{\phi} 
         & 
         \;H  
         \\
         FY \times GY \ar[ru]^{\lambda_{Y}} \ar@(r, dl)[rru]^{\phi_{Y}} 
        } 
$$
 such that $\phi
\circ \lambda_{X}\:=\: \phi_{X}$.
If $H$ is a locale, then there is a unique morphism of
locales $\cc{D}(D_{\Delta FG})^{\sim} \mr{\phi} H$ such that 
$\phi \circ \#\lambda_{X}\:=\: \phi_{X}$
\end{proposition}
\begin{proof}
1) follows immediately from 2) and 3).  2) and 3) by
construction of the poset $D_{\Delta FG}$ and
the fact that $\cc{D}(D_{\Delta FG})$ is the free inf-lattice on this
poset. The last statement follows then from lemma \ref{diaconescu}.
\end{proof} 

Notice that here (unlike in the case of functions between sets)
we do not abuse the notation and indicate the associated sheaf
morphisms \mbox{$\mathcal{D}(D_{\Delta FG}) \to lFunc(F, \,G)$} and 
$\mathcal{D}(D_{\Delta FG}) \to lBij(F, \,G)$ with the symbol $' \# '$. 

Given a
natural transformation  $F \mr{\sigma} G$, the corresponding point  
$lFunc(F, \,G) \mr{\sigma^{*}} 2$ is
characterized by: 
$$\sigma^{*}\# [(X,\,\langle \!x_{0}\,|\,x_{1}\!\rangle )] \;=\; 1 \iff 
                                                 \sigma X (x_{0}) = x_{1}$$

\vspace{2ex}

Next we prove the localic version of Yoneda's Lemma.

\begin{lemma} \label{yoneda}
Given any set valued functor $F: \cc{C} \mr{} \cc{S}$, and any object
$A \in \cc{C}$, functions  
$[A, \,X] \times FX \mr{\phi_{X}} FA$ defined by:
$$ 
Given \;\;  A\mr{x}X\,,  \;
y \in FX:\;\;  
\phi_{X}(x,\,y) = \{a \in FA \;|\; Fx(a) = y\}
$$ 
induce an isomorphism of locales $\phi: lFunc([A,\, -],\, F) \mr{\cong}
FA$ (where $FA$ denotes the discrete locale on the set $FA$).  
\end{lemma}
\begin{proof}
Given $X \mr{f} Y$, the equation   
$\phi_{X} \leq \phi_{Y} \circ (F(f) \times G(f))$ is immediate
to verify. It also becomes clear after inspection that the fact that  $Fx$
 is a function implies that the covering conditions 2.2 in proposition
\ref{genericnr} are satisfied. Thus, it follows there is a
(unique) morphism of  locales 
$\phi: lFunc([A,\, -],\, F) \mr{\cong} FA$ such that 
$\phi \circ \#\lambda_{X}\:=\: \phi_{X}$. We define a morphism of
locales $\lambda$ in the other direction by:
$$
Given \;\; I\subset FA:\;\;\;\; \lambda(I) = \bigvee_{a \in I} 
\#\lambda_A(id_A,\, a) = \bigvee_{a \in I}\# [(A,\,\langle \!id_A\,|\,a\!\rangle )]
$$
The reader can check that conditions 1.2.u) and 1.2.e) in
\ref{genericnr} imply that $\lambda$ preserves ``$\,\wedge\,$'' and
``$\,1\,$'' respectively. Since it clearly preserves ``$\,\bigvee\,$'', we
have that  $\lambda$ is a morphism of locales. Using this, we now
show that $\lambda$ is the inverse of $\phi$.

\vspace{1ex}

\emph{equation} $\phi \circ \lambda = id\,$): It is enough to show
for each $a \in FA$ that $\phi(\lambda\{a\}) = \{a\}$. But
$\phi(\lambda\{a\})  = \phi\#\lambda_A(id_A,\, a) = 
\phi_A(id_A,\, a) = \{a\} $, which is clear.

\vspace{1ex}

\emph{equation} $\lambda \circ \phi = id\,$): It is enough to show
for each $X \in \cc{C}$ the equation 
\mbox{$\lambda \circ \phi_X = \#\lambda_X$.} That is, for each
$A\mr{x}X\;and  \;y \in FX$: 
$\;\lambda \circ \phi_X(x,\, y) = \#\lambda_X(x,\, y)$:
$$
\lambda \circ \phi_X(x,\, y) \;=\; \lambda\{a \in FA \;|\; Fx(a) = y\} 
 \;=\;\;\bigvee_{a \in FA \;|\; Fx(a) = y}
\#[(A,\,\langle \!id_A\,|\,a\!\rangle )]\,.
$$
On the other hand, we have by \ref{genericnr} 1.2.e), 
$1 = \bigvee_{a \in FA}\#[(A,\,\langle \!id_A\,|\,a\!\rangle )]$. Taking infimun
against \mbox{$\#[(X,\,\langle \!x\,|\,y\!\rangle )]$} it follows 
$$
\#[(X,\,\langle \!x\,|\,y\!\rangle )] \;\;=\;\; \bigvee_{a \in FA}
                   \#[(A,\,\langle \!id_A\,|\,a\!\rangle )] \wedge
                   \#[(X,\,\langle \!x\,|\,y\!\rangle )]\,. \;\;\; 
$$ 
But   
$\#[(A,\,\langle \!id_A\,|\,a\!\rangle )] \;\leq\;
\#[(X,\,\langle \!x\,|\,Fx(a)\!\rangle )]$. Thus:
$$
\#[(A,\,\langle \!id_A\,|\,a\!\rangle )] \wedge \#[(X,\,\langle \!x\,|\,y\!\rangle )] \;=\;
               \#[(A,\,\langle \!id_A\,|\,a\!\rangle )] \;\;\;\;if \;Fx(a) = y
$$
$$
\#[(A,\,\langle \!id_A\,|\,a\!\rangle )] \wedge \#[(X,\,\langle \!x\,|\,y\!\rangle )] 
\;=\;\;\;\;\;\;\;\;\;\;\;\;0\;\;\;\;\;\;\;\;\;\;\;\;\;\;\;\;\; 
                                               \;\; if \;Fx(a) \neq y
$$
(the second equation using  proposition \ref{genericnr} 1.2.u). It follows
then 
$$
\#\lambda_X(x,\,y) \;=\; \#[(X,\,\langle \!x\,|\,y\!\rangle )]   \;\;=\;\;      
\bigvee_{a \in FA \;|\; Fx(a) = y}
\#[(A,\,\langle \!id_A\,|\,a\!\rangle )]\, .
$$
This finishes the proof of the equation $\lambda \circ \phi = id$.
\end{proof}

As usual, given any two objects $A, B \in \cc{C}$, it follows there is
an isomorphism of locales $
\phi: lFunc([A,\, -],\, [B, \,-]) \mr{\cong} [B, \,A]$. We also have:

\vspace{1ex}

\begin{lemma} \label{yonedaauto}
Given any category $\cc{C}$ and any two objects
$A, B \in \cc{C}$, the functions  
$[A, \,X] \times [B, \,X] \mr{\phi_{X}} Iso[B, \,A]$ defined by:
$$ 
Given \;\;  A\mr{x}X\,,  \;
A\mr{y}X:\;\;  
\phi_{X}(x,\,y) = \{B\mr{a}A \;\;|\;\; a\; iso \;\;and\;\; xa = y\}
$$ 
induce an isomorphism of locales $\phi: lBij([A,\, -],\, [B, \,-]) \mr{\cong}
Iso[B, \,A]$ (where $Iso[B, \,A]$ denotes the discrete locale on the set of
isomorphisms from $B$ to $A$).

In particular, we have an isomorphism of locales 
$\phi: lAut([A,\, -]) \mr{\cong} Aut(A)^{op}$ (where $Aut(A)$ denotes the
discrete locale on the set of automorphisms of $A$, and the $'op\,'$
indicates that there is a reversal of arrows in this last set).    
\end{lemma}
\begin{proof}
As in the proof of lemma \ref{yoneda}, the equation   
$\phi_{X} \leq \phi_{Y} \circ (F(f) \times G(f))$ is immediate
to verify. We leave the reader to inspect that the two
additional covering conditions 2.3 in proposition \ref{genericnr} follow readily
from the fact that the morphisms $B\mr{a}A$ are isomorphisms. Thus, it
follows there is a (unique) morphism of  locales 
$\phi: lBij([A,\, -],\, [B, \,-]) \mr{}
Iso[B, \,A]$
 such that 
$\phi \circ \#\lambda_{X}\:=\: \phi_{X}$. We define a morphism of
locales $\lambda$ in the other direction as in \ref{yoneda}: 
$$
Given \;\; I\subset Iso[B,\,A]:\;\;\;\; \lambda(I) = \bigvee_{a \in I} 
\#\lambda_A(id_A,\, a) = \bigvee_{a \in I}\# [(A,\,\langle \!id_A\,|\,a\!\rangle )]
$$
As in \ref{yoneda} it is straightforward to check that  $\lambda$ is a 
morphism of locales. Using this, we now show that $\lambda$ is the
inverse of $\phi$.

\vspace{1ex}

\emph{equation} $\phi \circ \lambda = id\,$): Same proof that in
\ref{yoneda}

\vspace{1ex}

\emph{equation} $\lambda \circ \phi = id\,$): As in \ref{yoneda} it is
enough to show, 
for each $X \in \cc{C}$, 
\mbox{$A\mr{x}X$} and $B\mr{y}X$, the equation
$\;\lambda \circ \phi_X(x,\, y) = \#\lambda_X(x,\, y)$:
$$
\lambda \circ \phi_X(x,\, y) \;=\; 
\lambda\{B\mr{a}A \;|\; a \; iso\,,\; xa = y\} \;=\;\;
\bigvee_{B\mr{a}A \;|\; a\; iso \,,\; xa = y}
\#[(A,\,\langle \!id_A\,|\,a\!\rangle )]\,.
$$
On the other hand, by the same reasoning as in \ref{yoneda} we have: 
$$
\#\lambda_X(x,\,y) \;=\; \#[(X,\,\langle \!x\,|\,y\!\rangle )]   \;\;=\;\;      
\bigvee_{B\mr{a}A \;|\; xa = y}
\#[(A,\,\langle \!id_A\,|\,a\!\rangle )]\, .
$$
Thus, to finish the proof we have to show that if $B\mr{a}A$ is not an
isomorphism, then $\#[(A,\,\langle \!id_A\,|\,a\!\rangle )] = 0$. We do this as
follows:

Notice that $B\mr{a}A$ is an isomorphism if and only if it is an
epimorphism and it has a left inverse $A\mr{x}B$, $xa = id_B$. Thus,
that $a$
is not an isomorphism means the following: 
$$
(\,\exists \;
\xymatrix{A \ar@<1ex>[r]^{y} \ar@<-1ex>[r]^{x} & X}  \;|\;\;\; 
x \neq y\,, \; xa = ya\,) \;\;\;\; or \;\;\;\;
(\,\forall \, A\mr{x}B\,, \; xa \neq id_B).\;
$$

Assume the first statement: By proposition \ref{yoneda} 1.3 i) it follows: 
$$\#[(A,\,\langle \!id_A\,|\,a\!\rangle )] \;\;\leq\;\; \#[(X,\,\langle \!x\,|\,xa\!\rangle )\,,\;
(X,\,\langle \!y\,|\,ya\!\rangle )] \;\;=\;\;  0\,.\;
$$ 

Assume the second statement:
By proposition \ref{yoneda} 1.3 s), 
$$\;\;\;\;\;\;\;\;\;\;\;1
\;\;=\;\;\bigvee_{A\mr{x}B}\#[(B,\,\langle \!x\,|\,id_B\!\rangle )]\;$$   
$$Thus \;\;\;\#[(A,\,\langle \!id_A\,|\,a\!\rangle )]\;\; = \;\; \bigvee_{A\mr{x}B}
\#[(B,\,\langle \!x\,|\,id_B\!\rangle )\,,\; (A,\,\langle \!id_A\,|\,a\!\rangle )] \,. $$ 
$$But \;\; \#[(B,\,\langle \!x\,|\,id_B\!\rangle )\,,\; (A,\,\langle \!id_A\,|\,a\!\rangle )] \;\leq\;  
\#[(B,\,\langle \!x\,|\,id_B\!\rangle )\,,\; (B,\,\langle \!x\,|\,xa\!\rangle )]\;$$
which is equal to $0$ for
all  $A\mr{x}B$ by proposition \ref{yoneda} 1.2 u). 
\end{proof}

\vspace{2ex}

Given any set valued functor $\cc{C} \mr{F} \cc{S}$, 
the localic space $lAut(F)$ is a localic group, and this group acts on
each set $FX$. 

More generally, given any two set valued functors 
$F,\;G\,:\cc{C} \mr{} \cc{S}$ and any object 
$X \in \mathcal{C}$, the map
$FX \times GX \mr{\lambda_{X}} \mathcal{D}(D_{F\Delta G})$
 determines morphisms of locales
$lRel(FX,\,GX) \mr{\lambda_{X}^{*}} lRel(F,\, G)$, 
$lFunc(FX,\,GX) \mr{\lambda_{X}^{*}} lFunc(F,\, G)$, and
\mbox{$lBij(FX,\,GX) \mr{\lambda_{X}^{*}} lBij(F,\, G)$,}  
defined by $\lambda_{X}^{*}[\langle \!x_{0}\,|\,x_{1}\!\rangle ] \;=\; 
\# [(X,\,\langle \!x_{0}\,|\,x_{1}\!\rangle )]$. These assertions follow since the map 
$FX \times GX \mr{\lambda_{X}} \mathcal{D}(D_{F\Delta G})$ sends
covers into covers on the respective sites of definition.

\vspace{1ex}

It is straightforward to check from
proposition \ref{groupstructure} the following:
 
\begin{proposition} \label{enrichement}
For any set valued functors $F,\; G,\; H\,:\cc{C} \mr{} \cc{S}$, there are
morphisms of locales:
$$
lFunc(F, \, G)\;\; \mr{m^{*}}\;\; lFunc(F, \, H) \otimes 
lFunc(H, \,G)\,,\ \;\;  lFunc(F, \, )\;\; \mr{e^{*}}\;\;2 
$$
defined on the generators by the following formulae:
$$m^{*}[(X, \;\langle \!x \:|\: y\:\rangle ] = (\lambda_{X}^{*} \otimes 
\lambda_{X}^{*})\, m_{X}^{*}[\langle \!x \:|\: y\:\rangle ]\,,\;\;\;
e^{*}[(X, \;\langle \!x \:|\: y\:\rangle ] =  \lambda_{X}^{*}e_{X}^{*}[\langle \!x \:|\:
y\:\rangle ]
$$
where $m_{X}^{*}$ and $e_{X}^{*}$ are the morphisms defined in
proposition \ref{groupstructure}.

These data satisfy the equations of an enrichment of the functor
category $\cc{S}^{\cc{C}}$ over the category
of localic spaces, which we denote $lFunc(\cc{S}^{\cc{C}})$. Clearly
the evaluation functor  $\cc{S}^{\cc{C}} \mr{ev_X}
\cc{S}$, $ev_X(F) = FX$ becomes a functor for the enriched
structures (by their very definition). This defines a (ordinary)
functor into the (ordinary)
category of enriched functors, which we denote $\mu^*$, $X \mapsto ev_X$,  
$\cc{C} \mr{\mu^*}  l\cc{S}^{lFunc(\cc{S}^\cc{C})}$.    


\vspace{1ex}

The above formulae together with 
\mbox{$\iota^{*}[(X, \;\langle \!x \:|\: y\:\rangle ] \;=\; 
                               \lambda_{X}^{*}\iota_{X}^{*}[\langle \!x \:|\: y\:\rangle ]$} 
define also morphisms of locales:
$$
lBij(F, \, G) \;\;\mr{m^{*}} \;\;lBij(F, \, H) \otimes lBij(H, \, G)
\,,\;\;\; lBij(F, \, F)\;\;\mr{e^{*}}\;\; 2
$$
$$
lBij(F, \, G)\;\; \mr{\iota^{*}} \;\;lBij(G, \, F)
$$
which determine a structure of localic groupoid on the (discrete) set
of objects of $\cc{S}^{\cc{C}}$, which we denote
$lBij(\cc{S}^{\cc{C}})$. As before, the evaluation functor becomes
enriched and defines an (ordinary) functor into the (ordinary)
category of enriched functors, which we also denote $\mu^*$, 
$X \mapsto ev_X$,  
$\cc{C} \mr{\mu^*}  l\cc{S}^{lBij( \cc{S}^\cc{C})}$. 
\qed \end{proposition}

Given a functor $\cc{C} \mr{T} \cc{H}$, any two set valued functors 
$F,\;G\,:\cc{H} \mr{} \cc{S}$, and any object 
$X \in \mathcal{C}$, there are morphisms of locales 
$lRel(FT,\,GT) \mr{T^*} lRel(F,\, G)$, 
$lFunc(FT,\,GT) \mr{T^*} lFunc(F,\, G)$, and
\mbox{$lBij(FT,\,GT) \mr{T^*} lBij(F,\, G)$,}
induced by morphisms of sites 
$\mathcal{D}(D_{FT\Delta GT}) \mr{T^*} \mathcal{D}(D_{F\Delta G})$
defined on the generators by the following formula: 
$T^*[(X, \;\langle \!x \:|\: y\:\rangle ] = [(TX, \;\langle \!x \:|\: y\:\rangle ]$ (that is, 
$T^*\lambda_{X} = \lambda_{TX}$). It is straightforward to check in
\ref{genericnr} that this map send covers into covers. Furthermore, it
is also straightforward to check that these data determine enriched functors
$lFunc(\cc{S}^\cc{H}) \mr{T^*} lFunc(\cc{S}^\cc{C})$, 
$lBij(\cc{S}^\cc{H}) \mr{T^*} lBij(\cc{S}^\cc{C})$,    
for the (co)-structure defined in \ref{enrichement}. Finally, it is
clear that all this is contravariantly functorial in the variable $\cc{C}$.
In all, this finishes the proof of (compare with \ref{transition}):

\begin{proposition} \label{transition1}
The assignments of the localic category $lFunc(\cc{S}^\cc{C})$ and of
the localic groupoid $lBij(\cc{S}^\cc{C})$ (with discrete set of
objects)  are contravariantly
functorial on $\cc{C}$, into the category of localic categories and
localic groupoids respectively. 

\qed \end{proposition} 

It follows immediately:

\begin{proposition}\label{transition2} 
The assignments of the categories of enriched functors together with
the functor $\mu^*$ in \ref{enrichement} are functorial in $\cc{C}$ in
such a way that $\mu^*$ becomes a natural transformation,  
$\cc{C} \mr{\mu^*}  l\cc{S}^{lFunc(\cc{S}^\cc{C})}$, 
$\cc{C} \mr{\mu^*}  l\cc{S}^{lBij( \cc{S}^\cc{C})}$. 
\qed \end{proposition}

\vspace{2ex}

\subsection{The localic groupoid of points of a topos}
\indent \vspace{1ex}
\label{lpoints}

Given any topos $\cc{E}$, consider a site $\cc{C}$ such that 
$\cc{E} = \cc{C}^{\sim}$. The usual
 equivalence between the category of
points of the site (that is, set valued flat and
continuous functors) and the category of (inverse images of) points
of the topos of sheaves $\cc{E} = \cc{C}^{\sim}$, induces an enriched structure
in the category and in the groupoid of points of $\cc{E}$. We define the 
\emph{localic groupoid of points} of a topos $\cc{E} = \cc{C}^{\sim}$ 
(which will be meaningful
in general only when the  topos has sufficiently many points)
to be 
$l\cc{P}oints(\cc{E}) \subset lBij(\cc{S}^{\cc{C}})^{op}$ (where
``$\subset$'' indicates the
enriched structure induced on the (full) subgroupoid whose objects are
the points of the site, changing the variance as in \cite{G2}). Notice
that given any two points $f, \,g$, we can define directly this
localic groupoid setting $l\cc{P}oints(\cc{E})[f,\,g] =
lBij(g^*,\,f^*)$. The given definition does not add rigor, but the
reason to do so is that it makes sense over an arbitrary base topos $\cc{S}$.
In the same way we define the \emph{localic category of points}.   

In this terminology and notation, from propositions \ref{transition1}
and \ref{transition2} it follows (recall that morphisms of topoi go on
the other way to the inverse image functors):
    
\begin{proposition} \label{comparison}
Let $\cc{E}$ be any topos. Then the assignment of the groupoid 
$l\cc{P}oints(\cc{E})$ is functorial in $\cc{E}$, into the category of localic
groupoids (with discrete set of objects) and morphisms of localic
groupoids. Furthermore, there is a geometric morphism of topoi $\mu$,
$\cc{B}(l\cc{P}oints(\cc{E})^{op}) \mr{\mu} \cc{E}$, whose inverse
image is given by $\mu^*(X) = ev_X$, and $\mu$ is natural in $\cc{E}$. 
\end{proposition}
\begin{proof}
The only point that needs some care is the existence of the geometric
morphism $\mu$. Consider any site $\cc{C}$, $\cc{E} = \cc{C}^{\sim}$. 
Clearly we have $\cc{C} \mr{\mu^*}
\cc{B}(l\cc{P}oints(\cc{E})^{op})$. Notice that
the family of points of
$\cc{B}(l\cc{P}oints(\cc{E})^{op})$ corresponding to evaluation
functors is surjective (\ref{betadiscreteobjects}). Then, it readily
follows that $\mu^*$ is flat and continuous.
\end{proof}

\vspace{2ex}

\subsection{The localic group of a pointed topos and filtered inverse limits} 
\indent \vspace{1ex}

We are interested in the functor $l\cc{P}oints$ behavior regarding
filtered inverse limits, but we left it for another
occasion (or the interested reader) the development of the general
theory (however, see comment \ref{comentario}). We now develop
in detail the particular case of pointed topoi and the
localic group of automorphism of the point. We also take care of some
necessary  \mbox{$2$-categorical} aspects, and for this purpose we first
review the results of the previous section in this particular case.  

\vspace{1ex} 

We have in particular for each point $\cc{S}
\mr{f} \cc{E}$, with $\cc{E} = \cc{C}^{\sim}$, 
$\cc{C} \mr{F} \cc{S}$, $F = f^{*}|_\cc{C}$,  the following:

\begin{proposition} \label{groupstructureF}
For any set valued
functor $F$, $\cc{C} \mr{F}  \cc{S}$, the localic space $lAut(F)$ is
a localic group which has
an action on the set $FX$ for each $X$, 
\mbox{$FX \times FX  \mr{\mu^{*}} lAut(F)$,} given by:
\mbox{$\mu^{*}(x_{0},\,x_{1}) \:=\: \#[(X,\,\langle \!x_{0}\,|\,x_{1}\!\rangle )]$},
and given any arrow $X \mr{f} Y$, the function  $FX \mr{F(f)} FY$ becomes a 
morphism of actions (see \cite{D}, 4.8).
This defines a lifting of $F$, which we denote $\mu F$, 
$\mathcal{C} \mr{\mu F} \mathcal{B}(lAutF)$. \qed   
\end{proposition}

\begin{proposition} \label{groupaction}
Let $F:\mathcal{C} \mr{} \cc{S}$ be any pointed site (that is,
$F$ is a flat and continuous functor inducing a point of the topos
 \mbox{$\cc{S} \mr{f} \cc{C}^{\sim})$.} Then, the 
lifting of $F$ defined in proposition \ref{groupstructureF} induces a
morphism of topoi, which we denote $\mu f$, commuting with the points:
$$
\diagrama
         {
          \cc{C}^{\sim} \ar[rr]^{\mu f^{*}} \ar[rd]^{f^{*}}  
          &&  \cc{B}(lAutF) \ar[ld]^{u{*}} 
          \\
          &  \cc{S}
         }
$$
\end{proposition}
\begin{proof}
Notice that
the canonical point of $\cc{B}(lAutF)$ is a surjection. Then, it readily
follows then that the lifting $\cc{C} \mr{\mu F} \cc{B}(lAutF)$ is flat and 
continuous.  
\end{proof}

\begin{proposition} \label{transition}
Given any two set valued functors related as in the diagram below:
$$
   \diagrama
     {
       \cc{C} \ar[rr]^{T} \ar[rd]^{F \;\;\;\;\; \theta} 
                                           &  &  \cc{D} \ar[ld]^{G}  \\
       & \cc{S} & & ,\; F \mr{\theta} GT \; \;(isomorphism).
     }
$$

1) There is a morphism of localic groups $lAut(G) \mr{lAut(T)} lAut(F)$
(induced by a morphism of sites as described below).

2) There is a morphism of sites (which we denote in the same way)  
\mbox{$\mathcal{D}(D_{\Delta F}) \mr{lAut(T)} \mathcal{D}(D_{\Delta G})$} 
defined by:
$$aut(T)\,[(X,\,\langle \!x_{0}\,|\, x_{1}\!\rangle )] 
\; =\;  [(TX,\,\langle \!\theta X(x_{0})\,|\,\theta X (x_{1})\!\rangle )]$$

3) If the functor $T$ has a left adjoint $S$, $id \mr{\eta} TS$, 
$ST \mr{\varepsilon} id$, then there is a natural transformation 
$G \mr{\sigma} FS$ which defines a morphism of inf-lattices 
$\mathcal{D}(D_{\Delta F}) \mr{aut(S)} \mathcal{D}(D_{\Delta G})$ 
by the formula;
$$aut(S)\,[(X,\,\langle \!x_{0}\,|\,x_{1})\!\rangle ] 
\; =\;  [(SX,\,\langle \!\sigma X(x_{0})\,|\,\sigma X
(x_{1})\!\rangle )]$$ Furthermore, $aut(S)$ is left adjoint to $aut(T)$.


\end{proposition}
\begin{proof}
1) The map between the localic spaces is given by the morphism of
sites defined in 2). It easily
follows from this definition and the definition of the (co)group structure 
(\ref{groupstructureF}, \ref{enrichement}) 
that the inverse image preserves this (co)structure.

2) We shall use \ref{genericnr}. 
Let $aut(T)^{*}_{X} : FX \times FX \mr{} \mathcal{D}(D_{\Delta G})$
   be defined by:
$$ 
 aut(T)_{X}(x_{0},\,x_{1}) 
\; =\; \lambda_{TX}((\theta X \times \theta X)(x_{0},\,x_{1}))   
    \; =\;  [(TX,\,\langle \!\theta X(x_{0})\,|\,\theta X (x_{1})\!\rangle )]
$$
Given $X \mr{f} Y$, by  naturality of $\theta$ it follows     
 $aut(T)_{X}\; \leq \; t_{Y} \circ (F(f) \times F(f))$. So, it
remains to check condition ii) in 2) of \ref{genericnr}. But this
immediately follows since $\theta X$ is a bijection (for the two basic
empty covers use injectivity, and for the two basic covers of $1$ use 
 surjectivity)

3) We shall use \ref{genericnr}. Let $G \mr{\sigma} FS$ be the composite 
$G \mr{G\eta} GTS \mr{\theta^{-1}} FS$. Then, as before, define
$aut(S)_{X} : GX \times GX \mr{} \mathcal{D}(D_{\Delta F})$ by:
$$ 
aut(S)_{X}(x_{0},\,x_{1}) 
\; =\; \lambda_{SX}((\sigma X \times \sigma X)(x_{0},\,x_{1}))   
    \; =\;  [(SX,\,\langle \!\sigma X(x_{0})\,|\,\sigma X (x_{1})\!\rangle )]
$$
It remains to prove that $aut(S)_{X}$ is left adjoint to
$aut(T)_{X}$.  This follows because
$\varepsilon$ and $\eta$ actually define arrows
in the poset $\cc{D}(D_{\Delta F})$. That this is the case amounts to
the validity of the equations
$F\varepsilon\,\circ\,\sigma T \,\circ \, \theta = id$, and 
$\theta S \circ \sigma = G\eta$. We verify this now:

\vspace{1ex}

$(F \mr{\theta} GT \mr{\sigma T} FST \mr{F\varepsilon} F)         \,=\,
 (F \mr{\theta} GT \mr{G\eta T} GTST \mr{\theta^{-1}ST} FST 
                                             \mr{F\varepsilon} F) \,=\,
 (F \mr{\theta} GT \mr{G\eta T} GTST \mr{GT\varepsilon} GT
                                               \mr{\theta^{-1}} F)  \,=\,
 (F \mr{\theta} GT \mr{\theta^{-1}} F) \,=\, (F \mr{id} F).$

The first equality by definition of $\sigma$, the second by
naturality, the third by the triangular equation of the adjointness,
and the fourth is obvious.
                 
\vspace{1ex}

$(G \mr{\sigma} FS \mr{\theta S} GTS) \,=\,
 (G \mr{G\eta} GTS \mr{\theta^{-1} S} FS \mr{\theta S} GTS) \,=\,
 (G \mr{G\eta} GTS).$

The first equality by definition of $\sigma$, and the second is
obvious.
\end{proof}

\begin{proposition} \label{mutheta}
In the situation of proposition \ref{transition}, assume $\cc{C}$ and
$\cc{D}$ are sites,  and $T$ a morphism of sites. Then, there is a
natural isomorphism 
\mbox{$\cc{B}(autT)^{*} \circ \mu f^{*} \mr{\mu \theta} \mu g^{*} \circ t^{*}$} 
as indicated in the following diagram:
$$        
\diagrama@1
        {
         \cc{C}\;\; \ar @<+2pt> `u[r] `[rr]^{\mu F} [rr]
                                             \ar[r]^{\epsilon} \ar[d]^{T}
         & \;\;\cc{C}^{\sim}\;\;  \ar[r]^{\mu f^{*}} 
                                  \ar[d]^{t^{*};\;\;\;\;\;\;\mu \theta}  
         & \;\;\cc{B}(lAutF) \ar[d]^{\cc{B}(autT)^{*}}
         \\
         \cc{D}\;\;  \ar @<-2pt> `d[r] `[rr]^{\mu G} [rr]  \ar[r]^{\epsilon} 
         & \;\;\cc{D}^{\sim}\;\;  \ar[r]^{\mu g^{*}}   
         & \;\;\cc{B}(lAutG)      
        }
$$
\end{proposition}

\vspace{1ex}

\begin{proof}
It is enough to define $\mu \theta$ as a natural transformation
\mbox{$\cc{B}(autT)^{*} \circ \mu F \mr{\mu \theta} \mu G \circ T$.} 
In order to do this, just check that given any object $X \in \cc{C}$, the
bijective function 
$FX \mr{\theta X} GTX$ is actually a morphism of actions.
\end{proof}

Triangles of set valued functors compose in the obvious way, and it is
straightforward to check that the constructions in the two propositions
above are functorial in the appropriate way. More precisely:

\begin{proposition}
Given two 
triangles and its composition:
$$
\diagrama
     {
      \cc{A}\;\; \ar[r] ^{R} \ar[rd]^{H} 
      & \;\;\cc{C}\;\; \ar[r]^{T} \ar[d]^{F} 
      & \;\;\cc{D} \ar[ld]^{G}
      && \cc{A}\;\; \ar[rr]^{S} \ar[rd]^{H} 
      &  
      & \;\;\cc{D} \ar[ld]^{G}   
      \\
      & \cc{S} 
      &&&& \cc{S}
     }
$$ 
fill respectively with natural isomorphisms $\xi$, $\theta$, and
$\kappa$ (where  
$S = T \circ R$ and  \mbox{$\kappa = \theta R \circ \xi$).} 

Then, 
$lAut(S)^{*} = lAut(T)^{*} \circ lAut(R)^{*}$, and
$\mu \kappa \,=\, \mu \theta r^{*} \circ \cc{B}lAut(T)^{*} \mu \xi$.    
\qed \end{proposition}

\vspace{2ex}

With this, we can state and prove the behavior regarding
filtered inverse limits.

\begin{proposition} \label{inverselimit}
Consider a filtered system of set valued functors and its colimit as
indicated in the diagram below: 
$$
\diagrama
 {
  \cc{C_{\alpha}} \ar[r]^{T_{\alpha \beta}} \ar[rd]^{F_{\alpha}}  
            &  \cc{C_{\beta}} \ar[d]^{F_{\beta}}  
            &  \cdots  \ar[r]  &  \cc{C}  \ar[lld]^{F} \\
       & \cc{S} 
     }
$$
Assume the triangles fill with natural isomorphisms 
$\theta_{\alpha \beta} : F_{\alpha} \mr{} F_{\beta}  T_{\alpha \beta}$
subject to the compatibility conditions 
$\theta_{\alpha \gamma} = \theta_{\beta \gamma}T_{\alpha \beta} 
                                     \circ \theta_{\alpha \beta}$
(it follows there are also natural isomorphisms 
$\theta_{\alpha} : F_{\alpha} \mr{} FT_{\alpha}$, where 
$\cc{C_{\alpha}} \mr{T_{\alpha}} \cc{C}$ are the inclusions into the
colimit). Then:

\vspace{1ex}

1)The induced (by \ref{transition}. 1)) cofiltered system of localic
groups: 
$$
  lAut(F_{\alpha}) \ml{t_{\alpha \beta}} lAut(F_{\beta})  
                                    \;\;\;\; \cdots \ml{}  lAut(F).
$$
is a cofiltered inverse limit of localic groups.

2) The induced (by \ref{transition}. 2) filtered system of
inf-lattices and site morphisms:
$$
\cc{D}(D_{\Delta F_{\alpha}}) \mr{t_{\alpha \beta}}
                          \cc{D}(D_{\Delta F_{\beta}}) \;\;\;\;
                                 \cdots  \mr{} \cc \cc{D}(D_{\Delta F}).
$$
is a filtered colimit of inf lattices, and the topology in 
$\cc{D}(D_{\Delta F})$ is the coarsest that makes the arrows 
$\cc{D}(D_{\Delta F_{\alpha}}) \mr{t_{\alpha}} \cc{D}(D_{\Delta F})$
 continuous. 
\end{proposition}
\begin{proof}
1) Follows immediately from 2) by lemma \ref{diaconescu}.

2) Here is where the filtering condition is necessary. An object of 
$\cc{C}$ is a germ of objects. That is, it is a pair $(X, \alpha)$,
with $X \in \cc{C_{\alpha}}$, two such pairs being considered equal if
they become equal further on the system. An arrow between two germs
is an arrow at some point in the system, two such arrows being
considered equal if they become equal further on the system. 
From this it readily follows that the objects of the inf-lattice 
$\cc{D}(D_{\Delta F})$ are germs of objects, and that the order
relation is what it should be. This shows that $\cc{D}(D_{\Delta F})$
is the filter colimit of the inf-lattices $\cc{D}(D_{\Delta
F_{\alpha}})$. It is immediate that the covers that generate the
topology in $\cc{D}(D_{\Delta F})$ are just the ones that generate the
coarsest topology which makes the arrows $t_{\alpha}$ continuous.  
\end{proof}

From Theorems \ref{inverselimittopoi} and \ref{inverselimitgroups} 
follows:  
\begin{proposition} \label{bottomrows}
In the situation of proposition \ref{inverselimit}, assume that each
$\cc{C}_{\alpha}$ is a site, each $T_{\alpha}$ a morphism of
sites, and $\cc{C}$ the inverse limit site (cf \ref{inverselimittopoi}).  
Assume furthermore that the transition morphisms $lAut(T_{\alpha \beta})$
given by proposition \ref{transition} are surjections. Then, in the
following diagram the two bottom rows are inverse limit diagrams of
topoi (where we use also the notation in propositions \ref{groupaction}
and \ref{mutheta}).
$$
\diagrama
      {
          \cc{C_{\alpha}}\;\; \ar[r]^{T_{\alpha \beta}} \ar[d]^{\epsilon}  
       &  \;\;\cc{C_{\beta}} \ar[d]^{\epsilon}  
       &  \cdots  \ar[r]  
       &  \;\;\cc{C}  \ar[d]^{\epsilon} 
       &  \cc{C_{\alpha}}\;\; \ar[r]^{T_{\alpha}} \ar[d]^{\epsilon} 
       &  \;\;\cc{C}  \ar[d]^{\epsilon} 
       \\
          \cc{C_{\alpha}^{\sim}}\;\;  \ar[r]^{t_{\alpha \beta}^{*}} 
          \ar[d]^{f_{\alpha}^{*}} \ar@{}[dr] |{\mu \theta_{\alpha \beta}}
       &  \;\;\cc{C_{\beta}^{\sim}} \ar[d]^{f_{\beta}^{*}} 
       &  \cdots  \ar[r]  
       &  \;\;\cc{C^{\sim}} \ar[d]^{f^{*}}
       &  \cc{C_{\alpha}^{\sim}}\;\; \ar[r]^{t_{\alpha}^{*}}
          \ar[d]^{f_{\alpha}^{*}} \ar@{}[dr] |{\mu \theta_{\alpha}} 
       &  \;\;\cc{C^{\sim}} \ar[d]^{f^{*}}
       \\
          \cc{B}lAut(F_{\alpha})\;\; \ar[r]^{lAut(T_{\alpha \beta})^{*}}  
       &  \;\;\cc{B}lAut(F_{\beta})   
       &  \cdots  \ar[r]  
       &  \;\;\cc{B}lAut(F)  
       &  \cc{C_{\alpha}}\;\; \ar[r]^{lAut(T_{\alpha})^{*}} 
       &  \;\;\cc{C} 
      }
$$
\qed \end {proposition}


\section{The fundamental theorems of galois theory} \label{galoistheory} 


The fundamental theorems of galois theory are representation theorems
for certain types of atomic topoi. We distinguish three cases in this
paper: the discrete case, corresponding to the classical galois
theory, the prodiscrete case, corresponding to Grothendieck's galois
theory, and the general localic case, that we call localic galois
theory.

\vspace{2ex}

\subsection{Pointed connected atomic sites} \indent \vspace{1ex}
\label{atomicsitesection} 

From the characterization of atomic sites given in \cite{BD} it is
easy to check the following:

\begin{proposition} \label{atomicsite}
Let $\mathcal{E}$ be a topos with a point 
 $\cc{S} \mr{f} \mathcal{E}$, and $\mathcal{C} \subset
 \mathcal{E}$ be a (small) full subcategory such that together with the 
 canonical topology is a pointed site $\mathcal{C} \mr{F}
 \cc{S}$ for $\mathcal{E}$, $\mathcal{E} = \mathcal{C}^{\sim}$,
$F = f^{*}|_{\cc{C}}$. Then:

\vspace{1ex}

1) If $\mathcal{E}$ is a pointed connected atomic topos, a
site as above can be chosen so that:

   i) Every arrow $Y \mr{} X$ in $\mathcal{C}$ is an strict 
epimorphism.

  ii) For every $X \in \mathcal{C}$ $FX \neq \emptyset$.

 iii) $F$ preserves strict epimorphisms.
 
  iv) The diagram of $F$, $\Gamma_{F}$, is a cofiltered category.

\vspace{1ex}

2) Given any pointed site as in 1), the topos of sheaves is a pointed
   connected atomic topos.

\vspace{1ex}

Condition ii) is equivalent to the connectedness of $\mathcal{E}$. The
category $\mathcal{C}$ can be taken to be the full subcategory of
non-empty connected objects, but not necessarily so.     

\qed \end{proposition}

The following two propositions are easy to prove (see \cite{D}).

\begin{proposition}  \label{poset}
The natural transformations $[Z, \:-] \mr{z^*} F$ are all injective,
and the diagram of $F$, $\Gamma_{F}$, is a cofiltered poset.
\qed \end{proposition}

\begin{proposition} \label{lfiel}
The functor $F$ is faithful (and reflects isomorphisms).
\qed \end{proposition}

\vspace{2ex}

\subsection{Discrete galois theory} \indent \vspace{1ex}
\label{Galois' galois theory} 

Discrete galois theory corresponds exactly to Artin's interpretation
of the classical galois
theory of roots of a polynomial with coefficients in a
field. We call this theory \emph{Galois' galois theory}, and its
fundamental theorem can be proved by elementary category methods (see
\cite{D}).
 The topos theoretical setting of this theory corresponds to the situation 
described in \ref{atomicsite} when the diagram $\Gamma_{F}$ of
the functor $F$ has a (co) final (i.e. initial) object, or,
equivalently, the inverse image functor of the point is 
representable. This means
(see \ref{yonedaauto}) that the localic group $lAutF$ is isomorphic to
the discrete group  
$Aut(A)^{op}$, where $A \in \cc{C}$ is any representing object. In
this case, the object $A$ is a \emph{universal covering} and the topos $\cc{E}$ in \ref{atomicsite} is said to be 
\emph{locally simple connected} (see \cite{BDLSCT}, where this notion was first
investigated in detail in the topos setting). Notice that since $A$ is
in $\cc{C}$, it is a cover, that is  $A \to 1$ is an epimorphism (this
is characteristic of the connected situation).

\begin{proposition} \label{basica} 
If $A$ is a representing object of $F$, every arrow 
$X \mr{f} A$ is an isomorphism. In particular, every endomorphism of 
$A$ is an isomorphism, $Aut(A) = [A, A]$, and if $G = [B,\, -]$ is any
other representable point, $A \cong B$.
\end{proposition}
\begin{proof}
By \ref{atomicsite} i) and iii),
it follows that there is $A \mr{g} X$
such that $fg = id$. Then, $g$ is a monomorphism. Since by
\ref{atomicsite} i) it is 
also a strict epimorphism, it follows that it is an isomorphism, and 
consequently so is $f$.
\end{proof} 
Let $\vartheta: [A,\,-] \mr{\cong} F$, with $A \in \cc{C}$ be a
representation of $F$, and let \mbox{$a = \vartheta A(id_A) \in FA$.} The
pair $(A,\;a)$ is an initial object in the diagram of $F$, and given
any $x \in [A,\, X]$, $\vartheta X(x) = F(x)(a)$. We have:


\begin{proposition}  \label{Asplit}
The object $A$ is a Galois object (see section \ref{GOsection}) and
every object $X$ is $A$-split with the set $[A, \, X]$. We have
$\cc{E} =  Split(A) \ml{\cong} \cc{P}_A$ in such a way that the given
point of $\cc{E}$ corresponds by the representing isomorphism
$\vartheta$ with the
canonical point of $\cc{P}_A$.  
\end{proposition}
\begin{proof}
Since $A$ represents a point, it is connected. 
Notice now that it is
enough to prove the statement for connected objects $X$. Let 
\mbox{$\theta:\;  \gamma^{*}[A,\, X] \times A  \mr{}  X \times A$} be the
arrow which corresponds under the adjunctions $\;\gamma^{*} \dashv
\gamma_{*}\;$ and $\;(-)\times A \dashv (-)^A\;$ to the arrow
$\;[A,\, X] \mr{} [A,\, X \times A]$, defined by $x \mapsto (x,
id_A)$. It can be seen that 
$F(\theta):\;  [A,\, X] \times FA  \mr{}  FX \times FA$ is given
by $(x,\, y) \mapsto (F(x)(a),\, y)$, that is $F(\theta) = (\vartheta
X,\, id_{FA})$. Thus, $F(\theta)$ is a bijection, and the proof
finishes by \ref{lfiel} (recall that by assumption \ref{atomicsite} i)
we already know that $A \to 1$ is a cover).
\end{proof}
Notice that since $A$ represents a point, it is not only a Galois object,
(thus a connected covering), but it is also projective, which means
that it is a universal covering.

\vspace{1ex}
 
In this representable case the fundamental theorems of galois
theory can be easily established. Clearly every set $[A,\, X]$
has an action of the group $G = Aut(A)^{op}$, thus the functor $F$
lifts into the topos of $G$-sets 
$\mathcal{C} \mr{\mu F} \mathcal{B}G$. 
It is not difficult to prove the following (see \cite{D} section 1):

\begin{theorem} \label{galois} For every object $X \in \mathcal{C}$
the action of the group $Aut(A)^{op}$ on the set $[A,\, X]$ is
transitive, and every arrow  $A \mr{x} X$  in $\mathcal{C}$ is the
categorical quotient by the action of the subgroup 
\mbox{$\{ h \in Aut(A) \: | \: xh = x \} \, \subset \, Aut(A)$} on $A$. 
\qed \end{theorem}

From this theorem, by easy general categorical arguments, follows:   
    
\begin{theorem}[\textbf{fundamental theorem}] \label{teorema0}
Let $\mathcal{E}$ be a pointed connected atomic topos 
 $\cc{S} \mr{f} \mathcal{E}$, and $\mathcal{C} \subset
 \mathcal{E}$ be a pointed site $\mathcal{C} \mr{F}
 \cc{S}$ for $\mathcal{E}$, $\mathcal{E} = \mathcal{C}^{\sim}$,
$F = f^{*}|_{\cc{C}}$ (in the sense of \ref{atomicsite} above) such
 that the functor $F$ is representable by an object $A \in \cc{C}$. Then, the
 lifting  $\mu F$ of $F$ lands in the subcategory of
 transitive $G$-sets, 
$\mathcal{C} \mr{\mu F} t\mathcal{B}G$, for the discrete group $G =
 Aut(A)^{op}$, and the induced morphism of
 topoi \mbox{$\cc{B}G \mr{\mu f} \cc{E}$}   is an equivalence.
\qed \end{theorem}

From this theorem, or by the same elementary proof, the
following groupoid version follows (compare with \cite{G2} Ex IV 7.6 d):

\begin{theorem}[\textbf{fundamental theorem}] 
\label{teorema00}
Let $\mathcal{E}$ be an essentially-pointed connected atomic topos,
and let
 $\cc{G}$ be its category of points (which is a connected
 discrete groupoid \ref{GG_Upointsunique}), $\cc{G} = \cc{P}oints(\cc{E})$. 
Then, the canonical
geometric morphism  \mbox{$ \cc{B}(\cc{G}^{op}) = \cc{S}^{\cc{G}^{op}}
\mr{\mu} \cc{E}$} is an equivalence.  
\qed \end{theorem}
 
We use the variance convention of SGA4. Given any geometric morphism
\mbox{$\cc{E} \mr{f} \cc{F}$,} clearly the induced morphism  
$\cc{B}(\cc{P}oints(\cc{E})^{op})  \mr{} \cc{B}(\cc{P}oints(\cc{F})^{op})$ 
makes
the square which expresses the naturality of $\mu$ 
commutative (here we have all ordinary categories, compare \ref{transition2}). 

\vspace{1ex}
This situation is characterized in terms of exactness
properties of the inverse image of the point. It is 
equivalent to the preservation of all limits by the inverse image
functor $f^{*}$, or, equivalently, the point is
essential (in the sense of \cite{G2}).
For simplicity we shall assume that $\cc{C}$ is the
full subcategory of all non empty connected objects.

\begin{proposition} [compare with \cite{G2} IV 7.6] \label{galois1}
In the situation of \ref{atomicsite}, assume that the functor 
\mbox{$\cc{E} \mr{f^{*}} \cc{S}$} preserves all (small) limits (the
point is essential). Then, the
diagram $\Gamma_{F}$ of the functor $\cc{C} \mr{F} \cc{S}$ has an
initial object $(A,\;a)$, and $F$ is representable by $A$.   
\end{proposition}  
\begin{proof}
Let B be the limit $B = Lim_{(X,\;x) \in \Gamma_{F}}  X$ taken in
$\cc{E}$. By assumption, the canonical morphism 
\mbox{$FB  \mr{}  Lim_{(X,\;x) \in \Gamma_{F}}  FX$} is a
bijection. Let \mbox{$a \in FB$} be the unique element corresponding under
this bijection to the tuple 
\mbox{$(x)_{(X,\;x) \in \Gamma_{F}} \in  Lim_{(X,\;x) \in \Gamma_{F}}
\; FX$}, and let A be the connected component of B so that $a \in
FA$. The verification of the statement in the theorem is standard
and straightforward.  
\end{proof}
\begin{corollary} \label{essentialiiflsc}
A pointed connected atomic topos is a locally simply connected Galois topos if
and only if the point is essential.
\qed \end {corollary}

\vspace{2ex}

\subsection{Prodiscrete galois theory} \indent \vspace{1ex} 
\label{Grothendieck's galois theory}

Grothendieck's galois theory corresponds to the situation
described in \ref{atomicsite} when the Galois objects are (co) cofinal in the
diagram $\Gamma_{F}$ of the functor $F$.
Then, by means of inverse limit techniques the fundamental theorem can be
proved by reducing it to the representable (or discrete) case. This yields a
prodiscrete localic group as the localic group of automorphisms of the
point. This is the method
introduced and  developed by Grothendieck in SGA1 \cite{G1} to treat the
profinite case (see also \cite{J}, and \cite{DC} for a detailed and
elementary description of all this). Later, in a series of commented
exercises in SGA4 \cite{G2} he gave
guidelines to treat the general prodiscrete case by means of locally
constant sheaves and progroups. The key result in these developments
is the construction of the Galois closure. 
In \cite{M2} I. Moerdijk developed this
program using prodiscrete localic groups
instead of progroups, and gave a rather sketchy proof of  the
fundamental theorem \mbox{(theorem 3.2 loc. cit.).} Prodiscrete 
localic groups and their classifying 
topoi are  completely equivalent to strict (in the sense that
transition morphisms are surjective) progroups and their classifying
topoi, as it was first observed by M. Tierney in lectures at Columbia 
University, and later stated in print independently by Moerdijk in
\cite{M2}. This result has been generalized to groupoids by Kennison
\cite{K}, 4.18 (see section \ref{locinvsection}).  

\vspace{1ex} 

Pointed  Galois topoi are given by pointed atomic sites (as explicitly
described in \ref{atomicsite}) so that pairs $(A, \:a)$, $a \in FA$, 
with $A$ a Galois object, are (co) cofinal in the diagram $\Gamma_{F}$ of $F$.

Let $(A,\;a)$ be an object in
the diagram of $F$, with $A$ a Galois object. Let $\cc{C}_{A}$ be the
full subcategory of $\cc{C}$ defined by:
$$
X \in \cc{C}_{A} \;\; \iff \;\;  [A,\;X] \mr{a^*} FX, \;  
a^*(h) = Fh(a), \;  is \;a \; bijection.
$$

\begin{proposition} \label{abajoA}  
An object $X \in \cc{C}$ is in $\cc{C}_{A}$ if and only if there exist an
arrow  $A \to X$.   
\end{proposition}
\begin{proof}
Consider an arrow $A \mr{x} X$. We have to see that $X \in
\cc{C}_{A}$. We have the commutative diagram:
$$
\diagrama
         {
          [A, \:A] \ar[r]^{a^{*}}   \ar[d]^{x_{*}}
          & FA \ar[d]^{F(x)} 
          \\
          [A, \:X] \ar[r]^{a^{*}}
          & FX
         }
$$
The bottom row is a bijection since it is already injective
(\ref{poset}) and $F(x)$ is surjective.
\end{proof}

Notice that it follows that $\cc{C}$ is the filtered  union (indexed
by the Galois objects in $\Gamma_{F}$) of the
full subcategories $\cc{C}_{A}$.

\vspace{1ex}

By definition $A \in \cc{C}_{A}$, and the restriction of the functor $F$ to  
 $\cc{C}_{A}$ is naturally isomorphic to $[A, \:-]$. Theorem
\ref{galois} gives:

\begin{proposition} \label{prelimit} 
The pair $\cc{C}_{A}$, $[A, \,-]$ defines an
atomic site with a representable point. The induced morphism 
$\cc{B}(lAutA^{op}) \mr{} \cc{E}_{A}$ is an equivalence (where $\cc{E}_{A}$
denotes the topos of sheaves on $\cc{C}_{A}$). 
\qed \end{proposition}

\vspace{1ex}


\begin{proposition} \label{galoisurjection}
Given a morphism $(B,\;b) \mr{f} (A,\;a)$ in $\Gamma_{F}$ with $A$
and $B$ Galois objects, there is a (full) inclusion of categories 
$\cc{C}_{A} \: \subset \: \cc{C}_{B}$. This gives rise to a triangle
as follows:
$$
\diagrama
         {
          \cc{C}_{A} \ar@{(->}[rr]  \ar[rd]^{[A, -]}  
           &&  \cc{C}_B \ar[ld]^{[B, -]}  
           \\
           & \cc{S} 
           && {f^*}: [A, \:-] \mr{\cong}  [B, \:-].
         }
$$
\end{proposition}
\begin{proof}
Let $X \in \cc{C}_{A}$. We have the commutative diagram:
$$
\diagrama
         {
          [A, \:-] \ar[rr]^{f^{*}} \ar[rd]^{a^{*}}
           && [B, \:-] \ar[ld]^{b{*}}
           \\
           & FX
         }
$$
Arrow $b^{*}$ is a bijection since it is already injective
(\ref{poset}) 
and, by assumption, $a^{*}$ is a bijection. Thus $X \in \cc{C}_{B}$. Clearly it
follows that $f^{*}$ is also a bijection. 
\end{proof}          

In lemma \ref{yonedaauto} we established an isomorphism 
$lAut([A, -]) \mr{\cong} Aut(A)^{op}$. We shall explicitly describe
 now how the morphism 
$lAut([B, -]) \mr{} lAut([A, -])$ defined in proposition \ref{transition}
is induced by a morphism $Aut(B) \mr{\varphi}  Aut(A)$. 

Let $l \in Aut(B)$, and let $h \in Aut(A)$ be the unique morphism such
that $f^{*} \circ h^{*} = l^{*} \circ f^{*}$ (recall that $f^{*}$ is an 
isomorphism). Define  $\varphi (l) = h$. Then, 
$\varphi (l) \circ f = f \circ l$ (since $f$ is an epimorphism, 
$\varphi (l)$ is characterized by this equation). Let $[(X,\;\langle \!x\,|\,y\!\rangle )]$
(with $A \mr{x} X$, \mbox{$A \mr{y} X$)} be a generator of the locale 
$lAut([A, \:-])$. Under the isomorphism with $Aut(A)$ it corresponds to
the open set $\{h \;|\; xh = y\}$ (where \mbox{$A \mr{h} A$).} Then, 
$\varphi^{-1}\{h \;|\; xh = y\} = \{l \;|\; x\varphi (l) = y\}$. Since $f$ is
an epimorphism, this set is equal to 
$\{l \;|\; x\varphi (l)f = yf\} =  \{l \;|\; xfl = y\} = \{l \;|\; f^{*}(x)l =
y\}$, which corresponds to $[(X,\;\langle \!f^{*}x\,|\,f^{*}y\!\rangle )]$. This
shows that the morphism of proposition \ref{transition} corresponds to 
$\varphi$ as defined above. Now, from $\varphi (l) \circ f = f \circ l
$ it follows \mbox{$F(\varphi (l)) \circ Ff = Ff \circ Fl$.} Since
$Ff(b) = a$     
we have \mbox{$F(\varphi (l))(a) = Ff \circ Fl(b)$,} that is   
$a^{*}(\varphi (l)) = Ff \circ b^{*}(b)$ (this equation also
characterizes $\varphi(l)$). Thus, 
\mbox{$\varphi = (a^{*})^{-1} \circ Ff \circ b^{*}$.}

We see, in particular, that $\varphi$ is then a surjective
function. This proves the following proposition:

\begin{proposition} \label{transitionsurjective}
The transition morphism between the localic groups corresponding to a
transition between two Galois objects in the diagram of $F$ is a surjection.

\qed \end{proposition}

\vspace{2ex}
 
We have the situation described in the following diagram:
$$
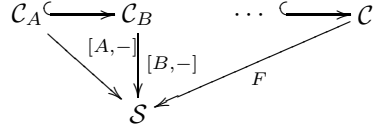
\diagrama
 {
  \cc{C}_{A} \ar@{^{(}->}[r]  \ar[rd]^{[A, -]}  
            &  \cc{C}_B \ar[d]^{[B, -]}  
            &  \cdots \;\; \ar@{^{(}->}[r]  &  \cc{C}  \ar[lld]^{F} \\
       & \cc{S} 
     }
$$
It follows from \ref{inverselimitgroups} that
the localic group $lAut(F)$ is prodiscrete and it is the inverse
limit of the induced filtered system of discrete groups $Aut(A)^{op}$
$$
  Aut(A)^{op} \ml{} Aut(B)^{op}  
                                    \;\;\;\; \cdots \ml{}  lAut(F).
$$ 

Furthermore, since $\cc{C}$ is the filtered  union (indexed by the Galois
objects in $\Gamma_{F}$) of the
full subcategories $\cc{C}_{A}$, and all the topologies are the
canonical one, $\cc{C}$ is the inverse limit site. It follows from 
\ref{galois} that the lifting 
$\cc{C} \mr{\mu F} \cc{B} lAut(F)$ of $F$ 
(\ref{groupaction}) lands in the subcategory of
 transitive $G$-sets. Furthermore, from  \ref{inverselimittopoi}, 
 \ref{transitionsurjective} and \ref{bottomrows} it follows that both rows in
the following diagram are filtered inverse limits
of topoi (indexed by the Galois objects in $\Gamma_{F}$):

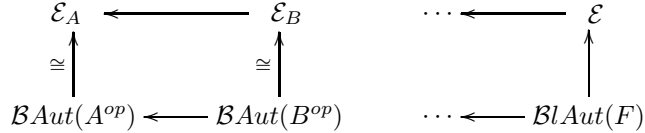
$$
\diagrama
      {
          \cc{E}_{A}\;\;    
       &  \;\;\cc{E}_{B} \ar[l]   
       &  \cdots    
       &  \;\;\cc{E}  \ar[l]  
       \\
          \cc{B}Aut(A^{op}) \ar[u]^{\cong}   
       &  \cc{B}Aut(B^{op}) \ar[l] \ar[u]^{\cong}   
       &  \cdots    
       &  \cc{B}lAut(F) \ar[l] \ar[u]  
      }
$$

Therefore the arrow $\cc{B} lAut(F) \mr{} \cc{E}$ is also an equivalence.
In conclusion, this finishes the proof of  the following theorem:

\begin{theorem}[\textbf{fundamental theorem}] \label{thSGA4}
Let $\mathcal{E}$ be a pointed Galois topos 
 \mbox{$\cc{S} \mr{f} \mathcal{E}$,} and $\mathcal{C} \subset
 \mathcal{E}$ be a pointed site $\mathcal{C} \mr{F}
 \cc{S}$ for $\mathcal{E}$, $\mathcal{E} = \mathcal{C}^{\sim}$,
\mbox{$F = f^{*}|_{\cc{C}}$} (in the sense of \ref{atomicsite}
 above). Then, the 
localic group $G = lAut(F)$ is prodiscrete, the lifting 
$\cc{C} \mr{\mu F} \cc{B} lAut(F)$ of $F$ 
(\ref{groupaction}) lands in the subcategory of
 transitive $G$-sets, and the induced morphism of
 topoi \mbox{$\cc{B}lAut(F) \mr{\mu f} \cc{E}$}   is an equivalence.
\qed \end{theorem}
Recall that if $\cc{S} \mr{f} \cc{E}$ is the corresponding point of
 the topos, $F = f^{*}|_{\cc{C}}$, then $lAut(F) = lAut(f)^{op}$.

\vspace{1ex}

From this theorem follows a groupoid version:

\begin{theorem}[\textbf{fundamental theorem}] 
\label{teorema0000}
Let $\mathcal{E}$ be a Galois topos with points (thus enough). Then the
canonical geometric morphism \mbox{$\cc{B}(\cc{G}^{op}) \mr{} \cc{E}$} is an
equivalence, where $\cc{G}$ is the localic groupoid of
points $\cc{G} = l\cc{P}oints(\cc{E})$ defined in section
\ref{lpoints}, and this groupoid has prodiscrete ``hom''spaces (in
particular, prodiscrete vertex localic groups). 
\end{theorem}
\begin{proof}
Let $\cc{S} \mr{f} \cc{E}$ be any point of $\cc{E}$, and consider
the commutative diagram:
$$
\xymatrix
         {
            \cc{B}(lAut(f)^{op}) \ar[rr] \ar[rd]
          & 
          & \cc{B}(\cc{G}^{op}) \ar[ld]
          \\
          & 
            \cc{E}  
         }
$$
From the second statement of 
\ref{gsetatomic1} and \ref{thSGA4} it follows that the
horizontal arrow and the left diagonal are equivalences. Thus the
remaining arrow is so  (compare with \cite{G1} V 5.8). The last
statement follows from \ref{comentario} below.         
\end{proof}  

The reader
should be aware that the groupoid in this theorem is not a prodiscrete
localic groupoid.

\vspace{1ex}

We now characterize this situation in terms of exactness
properties of the inverse image of the point. Theorem
\ref{cotensortheorem} below is inspired in a natural way of
constructing a normal 
covering (which covers a given covering) in the classical topological
theory of covering spaces. In fact, this theorem is an explicit
construction of the Galois closure.    

Grothendieck's theory corresponds to the case in which the
point, although not necessarily essential, is such that the inverse
image preserves certain infinite limits, namely, cotensors of
connected objects. This is equivalent to the existence of Galois
closure (that is, the Galois objects generate the topos), or to the
fact that the localic group $lAut(p)$ is prodiscrete. We 
elaborate on this now.

Consider a pointed connected atomic topos $\cc{S} \mr{f} \cc{E}$ and
a corresponding pointed site $\cc{C} \mr{F} \cc{S}$ as in
\ref{atomicsite}. For simplicity we shall assume that $\cc{C}$ is the
full subcategory of all non empty connected objects. Recall that the
topology is the canonical topology.

\begin{definition} \label{cotensoraxiom}
Let $\mathcal{E} \mr{\gamma} \cc{S}$ be any topos. We say that a
point $\cc{S} \mr{p} \mathcal{E}$ of $\cc{E}$ is 
\emph{proessential} if the
inverse image preserves cotensors of connected objects. That is, given
any connected object $X$ and any set $S$, the canonical morphism:
$$ p^{*}(X^{\gamma^{*}S}) = p^{*}(\prod\nolimits_{\,S} X)  \mr{}  
                               \prod\nolimits_{\,S} p^{*}X = (p^{*}X)^{S}. $$
is a bijection.
\end{definition}

Notice that preservation of cotensors (of any object)
is a much stronger condition which implies that the point is
essential (see \cite{BP}). 

\begin{theorem} \label{cotensortheorem}
In the situation of \ref{atomicsite}, assume that the point 
$\cc{S} \mr{p} \cc{E}$, \mbox{$F = p^*|_{\cc{C}}$} is proessential. Then, 
the objects $(A,\;a)$ with $A$ a Galois object are cofinal in the diagram
$\Gamma_{F}$ of $F$. In fact, the following holds: Given any connected
object $X$, there exists a Galois object $A$, and an element $a\in FA$
such that for all $x \in FX$ there exists an arrow $(A,\;a) \mr{f}
(X,\;x)$ such that $Ff(a)= x$. Notice that in this context $f$ is
unique and a strict epimorphism (\ref{poset} and \ref{atomicsite} i).
\end{theorem}
\begin{proof}
Let B be the cotensor $B = \prod_{FX} X$ taken in $\cc{E}$. By
assumption, 
the canonical morphism 
\mbox{$F(\prod_{FX} X) \mr{} \prod_{FX} FX$} is a
bijection. Let $a \in FB$ be the unique element corresponding under
this bijection to the tuple 
\mbox{$(x)_{x \in FX} \in  \prod_{FX} FX$}, and let A be the
connected component of B such that $a \in FA$. Clearly, for each $x \in
FX$ there is an arrow 
in $\Gamma_{F}$ given by the projection $A \mr{\pi_x} X$,
characterized by the equation $F\pi_x(a) = x$. We prove now
that A, with the element $a \in FA$, is a Galois object. To this end we 
establish:

\vspace{1ex}

\emph{Lemma}. Given any $b \in FA$, there exists $A \mr{f_x} X$ such
that $Ff_x(b) = x$.

\vspace{1ex}

Clearly, from this it follows (by the universal property of the
product) that there exists  
$A \mr{h} A$ such that $Fh(b) = a$. Let $c = Fh(a)$, and apply the
lemma to this element $c \in FA$. It follows as before that there is
$A \mr{g} A$ such that $Fg(c) = a$. Then, by \ref{poset} it must be 
$g \circ h = id$. So h is a monomorphism, and thus by \ref{atomicsite}
i) it is an isomorphism. This shows that $A$ is a Galois object.

\vspace{1ex}

\emph{proof of the lemma}. 
Consider the action $\mu$ of $lAut(F)$ (\ref{groupaction}).
Take any $x \in FX$. Then:
$$\;\; 1 \;=\; \bigvee\nolimits_{z \in FX}\, \mu^{*}[\langle \!z \:|\: x\!\rangle ]$$
Since the action is transitive (\ref{transitivity2}), taking the
infimun against $\mu^{*}[\langle \!a \:|\: b\!\rangle ]$ yields:
$$ 0 \;\neq\; \mu^{*}[\langle \!a \:|\: b\!\rangle ] \;\leq\; 
\bigvee\nolimits_{z \in FX}\, \mu^{*}[\langle \!a \:|\: b\!\rangle ] \wedge 
                                             \mu^{*}[\langle \!z \:|\:x\!\rangle ]$$
It follows that there exists $z \in FX$ such that 
$$0 \,\neq\, \mu^{*}[\langle \!a \:|\: b\!\rangle ] \wedge \mu^{*}[\langle \!z\:|\:x\!\rangle ]$$
Since $\pi_z$ is a morphism of actions, we have:  
$$\mu^{*}[\langle \!a \:|\: b\!\rangle ] \;\leq\; \mu^{*}[\langle \!\pi_za \:|\: \pi_zb\!\rangle ] 
                                  \;=\; \mu^{*}[\langle \!z \:|\:\pi_zb\!\rangle ]$$
It follows:
$$0 \,\neq\, \mu^{*}[\langle \!z \:|\: \pi_zb\!\rangle ] \wedge \mu^{*}[\langle \!z\:|\:x\!\rangle ] 
     \;=\; \mu^{*}([\langle \!z \:|\: \pi_zb\!\rangle ] \wedge [\langle \!z\:|\:x\!\rangle ])$$
Thus, $0 \,\neq\,[\langle \!z \:|\: \pi_zb\!\rangle ] \wedge [\langle \!z\:|\:x\!\rangle ])$,
which implies $x = \pi_zb$. We set $f_x = \pi_z$. This finishes the
proof of the lemma.  
\end{proof}
 On the other hand, given any pointed Galois topos, it is easy to see
 that the point is proessential. In fact, given any
 connected object $X$, take a Galois object $A$ such that $X \in
 \cc{C}_A$. Any cotensor of $X$ lives in $\cc{E}_A =
 \cc{C}_{A}^{\sim}$, and the result follows since the restriction of
 the inverse image to the full subcategory $\cc{C}_A$ determines an
 essential point of $\cc{E}_A$. We have:  
\begin{corollary}
A pointed connected atomic topos is a Galois topos if and only if the point is
proessential.
\qed \end {corollary}

\vspace{2ex}

\subsection{Unpointed prodiscrete galois theory} 
\indent \vspace{1ex}\label{unpointedsection}

  
M. Bunge \cite{B} (see also
\cite{BM}) developed an unpointed theory for Galois topoi based on 
Joyal-Tierney descent theory \cite{JT}, and following Grothendieck's inverse
limit techniques along the lines of the pointed theory of \cite{M2},
necessarily in this case in terms of localic groupoids. Around the same time,
J. Kennison \cite{K} also developed an unpointed theory with a
different approach, but we shall not elaborate on this theory here.

We now describe briefly the unpointed theory along the lines of 
\cite{B}, \cite{BM} section 2, and describe explicitly the fundamental
groupoid as the localic groupoid of ``points'' (which may not be there !).  
We do this in a independent way of the representation theorems in \cite{JT}.
We shall show that Grothendieck's theory of sections
\ref{Galois' galois theory} and 
\ref{Grothendieck's galois theory} can as well be developed in a
pointless way. We also think that it will interest the reader to see
explicitly how the points, which are always there at the starting line
(the topoi $Split(U)$ always have points), are lost along the way.

\vspace{1ex}

Consider a Galois Topos as in definition \ref{galoistopos}. In section
\ref{Grothendieck's galois theory} the point furnish a filtered poset
$\Gamma_F$ along which compute an inverse 
limit of pointed topoi. In the absence of the point we have to deal
differently. Proposition  \ref{filtergalois} is at the base of this
development. Even though all the topoi in the system furnished by this
proposition have points, the
system is not a pointed system in the sense that there is no 
simultaneous choice of points commuting with the transition
morphisms. In fact, such a choice is equivalent to a point of the
inverse limit topos.

\vspace{1ex}

Consider now any connected locally connected topos $\cc{F}$ and the
Galois topos \mbox{$\cc{E} = GLC(\cc{F})$} (notice that $\cc{E}$ can be any
Galois topos \ref{GLCisgalois}).   
Given a morphism between Galois objects $A \mr{f} B$, the geometrical
morphism $Split(A) \ml{} Split(B)$ (with inverse image the full
inclusion of categories) clearly induces a surjective function between
the sets of points (for the surjectivity compare with
\ref{galoisurjection}) .  A point of
$\cc{E}$  furnish a way of choosing a point
(consistently with respect to the transition morphisms) on each topoi
$Split(A)$, thus, it is exactly an element of the inverse limit of
the sets of points of the topoi $Split(A)$. This inverse limit may be
empty, but taken in the category
of localic spaces it always defines a non trivial prodiscrete localic 
space (since the projections are surjective \cite{JT} IV 4.2.) $G_0$,
which is the space of (may be phantom) points
of the inverse limit topos $\cc{E}$. 

More over, there is induced a groupoid morphism
\mbox{$\cc{G}_A \ml{} \cc{G}_B$} between the categories (which are discrete
connected groupoids) of points, 
\mbox{$\cc{G}_A = \cc{P}oints(Split(A))$} (compare \ref{comparison}).
 The inverse limit of this filtered
diagram, taken in the category of localic groupoids, defines a 
 \emph{prodiscrete localic groupoid} (see \cite{K}, definition 2.8)
$\cc{G}$ with the prodiscrete localic space $G_0$ as its localic space of
`\mbox{`objects''.}  

We shall say that $\cc{G}$ is the \emph{localic groupoid of phantom
points} of the
Galois topos $\cc{E}$, and write $ph\cc{P}oints(\cc{E})$.  The points
(if any) of $G_0$ are exactly the points of the topos $\cc{E}$. 
On the other hand, there are also geometric morphisms between the
push-out topoi $\cc{P}_A \ml{} \cc{P}_B$, which are morphisms of
pointed topoi for the canonical points (cf \ref{U=A} and
\ref{GG_Upointsunique}). Thus, there is always a consistent choice of
points for the system of push-out topoi $\cc{P}_A$. 
      
The whole situation we have at hand is synthesized in the
following diagram:

$$
\xymatrix
        {
         \cc{G}_A   
        &\cc{G}_B  \ar[l] 
        & \;\; \cdots \ar[l] 
        &\cc{G} 
       \\
          \cc{B}(\cc{G}_{A}^{op})  \ar[d]  
        & \cc{B}(\cc{G}_{B}^{op})  \ar[l] \ar[d]
        & \;\; \cdots \ar[l] 
        & \cc{B}(\cc{G}^{op})  \ar[d]
       \\
          Split(A)  \ar[d]  
        & Split(B)  \ar[l] \ar[d]
        & \;\; \cdots \ar[l] 
        & \cc{E}  \ar[d]
       \\
          \cc{P}_A    
        & \cc{P}_B  \ar[l]  
        & \;\; \cdots \ar[l]
        & \cc{P}  
       \\   
        & \cc{S} \ar[lu] \ar[u] \ar[rru]
       }
$$
\vspace{1ex}
The isolated first row is an inverse limit by definition. That the
second row 
is an inverse limit means that
the functor $\cc{B}$ commutes with the inverse limit of discrete
groupoids which defines $\cc{G}$ in the first row. This is
proved in \cite{K} 4.18. That the third row is
an inverse limit is proposition \ref{filtergalois}. Finally, we 
define $\cc{P}$ as the mathematical object which makes the fourth row an
inverse limit. In the previous considerations we already saw that
everything commutes in the appropriate way, and this implies
the existence of the point $\cc{S} \mr{} \cc{P}$. 

The vertical down arrows on the left of the dots  are
equivalences by \ref{teorema00} and \ref{GG_U} respectively. It
follows (using  the horizontal rows) that the arrow 
\mbox{$\cc{B}(\cc{G}^{op}) \mr{} \cc{E}$} is an equivalence. This finishes the
proof of:

\begin{theorem}[\textbf{fundamental theorem}] 
\label{teorema000}
Let $\mathcal{F}$ be any connected locally connected topos, and
$\cc{E}$ be the Galois topos $\cc{E} = GLC(\cc{F})$. Then the
canonical geometric morphism \mbox{$\cc{B}(\cc{G}^{op}) \mr{} \cc{E}$} is an
equivalence, where $\cc{G}$ is the prodiscrete localic groupoid of
(phantom) points $\cc{G} = ph\cc{P}oints(\cc{E})$ defined
by the inverse limit above. 
\qed \end{theorem}  

\vspace{1ex}

Now, since each $Split(A)$ is equivalent to
$\cc{P}_A$, $\cc{E}$ should be equivalent to $\cc{P}$, and it
would follow then that the topos $\cc{E} = GLC(\cc{F})$ (and so any Galois
topos) always has a point ?. 

The problem here is that the system  of push
 out topoi is not a filtered system and \emph{can not be used as
 such} to define an inverse limit topos. Given $A \leq B$ in $GCov(\cc{E})$,
 the transition morphism $\cc{P}_A \ml{} \cc{P}_B$, which now we shall
 denote $p_f$, depends on the arrow  $A  \mr{f} B$ which witnesses that 
 $A \leq B$. The reader can easily verify this  by direct inspection. 
However, there is no
 complete chaos here. Given any two
 arrows $f,\; g\,: \, A \mr{} B$, it is also immediate to check by the
 same method that there is a 
(canonical) invertible natural transformation $p_f \,\cong\,
 p_g$, and that all these two-cells define a \emph{biordered
 inversely bifiltered two system} (see \cite{K} for this notion and
 further references) of topoi which is not inversely
 filtered. It is not known if the inverse limit (or bilimit) of such a
 thing is a topos, and even less what kind of topos if that were the
 case. The equivalences $Split(A) \mr{} \cc{P}_A$ induce an arrow on
 the inverse limits (whatever it is $\cc{P}$) 
 $\cc{E} \mr{} \cc{P}$ which presumably  will not have a pseudo
 inverse to compose with the point of $\cc{P}$ to give a point for 
 $\cc{E}$. 

\vspace{2ex}

\subsection{Comparison between the pointed and unpointed theories} 
\indent \vspace{1ex}\label{comparisonsection}

In the presence of points, both the pointed and the
unpointed theories apply, but do not furnish the same groupoid in the
fundamental theorem.  We now study how the two constructions of localic 
groupoids are related. Namely, the localic groupoid 
$l\cc{P}oints(\cc{E})$ in section 
\ref{Grothendieck's galois theory}, and the localic groupoid
$ph\cc{P}oints(\cc{E})$ in section \ref{unpointedsection}. 
It turns out that both correspond to filtered inverse limits of
discrete groupoids, but taken in different categories.   

This concerns the preservation 
of the filtered inverse limit of topoi 
$$
       \xymatrix
               {          
                  Split(A)    
                & Split(B)  \ar[l] 
                & \;\; \cdots \ar[l] 
                & \cc{E} 
               }
$$
by the functor $l\cc{P}oints$ defined in \ref{comparison}. As we shall
 see, this inverse limit is preserved into the category of localic
 groupoids with discrete space of objects. 

First, notice that since all the points of the essentially pointed
topoi $Split(A)$ are representable, it follows from the localic
Yoneda's lemma \ref{yonedaauto} that the localic and the discrete
groupoids 
of points are equivalent in this case. We
have \mbox{$\cc{G}_A = \cc{P}oints(Split(A)) \cong l\cc{P}oints(Split(A))$.}  
   
Given any two points  $f, \, g$ of $\cc{E}$, they are given by
compatible (with respect to the transition morphisms) tuples 
$f = (f_A)$, $g = (g_A)$ of points of the topoi $Split(A)$. Consider
the filtered system of discrete spaces (where $\cc{C}_A = Split(A)$).
$$
       \xymatrix
               {          
                  l\cc{P}oints(\cc{C}_A) [f_A, \,g_A]    
                & l\cc{P}oints(\cc{C}_B) [f_B, \,g_B]  \ar[l] 
                & \;\; \cdots \ar[l] 
                & l\cc{P}oints(\cc{E})[f,\,g] 
               }
$$
Taking into account the proof of \ref{comparison}, with the same
arguments as in the proof of proposition \ref{inverselimit} (where the case of
one of the vertex localic groups is donned in detail), it
follows that this diagram is an inverse limit diagram of localic
spaces.

This shows that $l\cc{P}oints(\cc{E})$ is the inverse limit of the
filtered system of discrete groupoids 
$\cc{G}_{A} \ml{} \cc{G}_{B} \; \cdots \;$ in the category of localic
groupoids with discrete space of objects, while
$ph\cc{P}oints(\cc{E})$ is by definition the inverse limit of the same
system in the category of all localic groupoids. It follows then that
there is a comparison morphism of localic groupoid 
$l\cc{P}oints(\cc{E}) \mr{}  ph\cc{P}oints(\cc{E})$. 

Notice that from
the representation theorems \ref{teorema000} and \ref{teorema0000}
follows that this morphism induces 
an equivalence between the classifying topoi.

\begin{comment} \label{comentario}
In the arguing above it is given an sketch of the
proof that the functor $l\cc{P}oints(\cc{E})$ preserves filtered
limits of topoi (into the category of localic groupoids with discrete
space of objects), generalizing \ref{inverselimit}. 
\end{comment}
\vspace{1ex}

Proposition \ref{betadiscreteobjects} says in a way that the
classifying topos of a localic groupoid with discrete space of objects
is a rather simple construction, similar to the classifying topos of 
a discrete groupoid. Based on this, it can be proved that when the inverse
limit topos $\cc{E}$ has points, the functor
$\cc{B}$ preserves the filtered inverse limit of the system of
localic groupoids $\cc{G}_A$ considered above (limit taken in the category
of localic groupoids with discrete space of objects). Use the fact that all the points
$f_A$ of the topoi $Split(A)$ are representable by projective objects,
and that this implies that
the transition morphisms  
$\cc{P}oints(\cc{C}_A) [f_A, \,g_A] \ml{} \cc{P}oints(\cc{C}_B) [f_B,\,g_B]$ 
are surjective (compare with \ref{inverselimitgroups}). It follows
that the topos
$\cc{B}l\cc{P}oints(\cc{E})$ is the inverse limit of the system of topoi 
$\cc{B}(\cc{G}_{A}^{op})$ (as it was the case for the topos 
$\cc{B}ph\cc{P}oints(\cc{E})$ by \cite{K} 4.18). This gives a proof of
\ref{teorema0000} along the lines of the proof of \ref{teorema000}, and
without the need of using \ref{gsetatomic1}. At the same time it shows
directly why the comparison morphism between the two localic groupoids
of points induces an equivalence of the classifying topoi.

\vspace{1ex}

 

Galois topos with points are connected but may have a non-connected
groupoid of points.
We finish this section with an example:

\begin{example} \label{example1} 
In SGA4 IV 7.2.6 d) it is said that there exists a strict progroup 
$H = (H_i)_{i \in I}$  such that the classifying topos $\cc{B}H$ has 
two non isomorphic points. Equivalently, there is a prodiscrete localic
group $H$ such that $\cc{B}H$ has two non isomorphic points.
This implies that the groupoid of points $\cc{P}oints(\cc{B}H)$ is not
connected. Its classifying 
topos $\cc{B}\cc{P}oints(\cc{B}H)$ can not be $\cc{B}H$. However, the
localic groupoid of points $l\cc{P}oints(\cc{B}H)$ (which has discrete
set of objects) is connected (in particular, the localic
space of morphisms between any two points is non trivial). We
have $\cc{B}H \cong \cc{B}l\cc{P}oints(\cc{B}H)$ and 
$H \cong l\cc{P}oints(\cc{B}H)$
\qed \end{example}

\vspace{2ex}
   
\subsection{Localic galois theory} \indent \vspace{1ex}
\label{localicgalois}

In the previous sections we have developed the fundamentals of the
galois theory as given by Grothendieck's guidelines up to its natural
end point, which is the representation theorems of Galois Topoi
\ref{thSGA4} and \ref{teorema000}. One aspect of
these theorems is that they furnish an axiomatic
characterization of the classifying topoi of prodiscrete
localic groups  
and (connected) prodiscrete localic groupoids respectively. 
With the notion of localic group generalizing the
notion of progroup, the natural end point of the theory is push
forward into the representation theorem of pointed connected atomic
topoi, which would be, in particular, an  
axiomatic characterization of the classifying topoi of general
localic groups. This theorem is \cite{JT} Ex.VIII 3. Theorem 1,
and it still generalizes closely Grothendieck's galois theory. In
particular, the localic group in the statement is still the localic
group of automorphisms of the point (or the localic groupoid of all points
as defined in section \ref{lpoints}), and as such, it is canonically
associated to the topos (and functorial).

\vspace{1ex}
   
 We now recall the fundamental theorems of \emph{localic galois
theory} established in \cite{D}, where the 
representation theorem of pointed connected atomic topoi is a
consequence of a 
theory completely different to the Joyal-Tierney theory, and more
akin in its methods to classical
galois theory (compare with section \ref{Galois' galois theory})

\vspace{1ex}

 Assume that the pair 
$\mathcal{C}$, $\mathcal{C}  \mr{F} \cc{S}$,
is a pointed connected atomic site in the sense explicitly described
in \ref{atomicsite} above. We have:

\begin{theorem} \label{transitivity2}
For every object $X \in \mathcal{C}$ the action of the localic group 
of automorphisms  $lAut(F)$ on the set $FX$ is transitive. That is,
given any pair  \mbox{$(x_{0},\,x_{1}) \in FX \times FX$}, 
$\#[(X,\,\langle \!x_{0}\,|\,x_{1}\!\rangle )] \; \not= \; 0$ in $lAut(F)$. 
\qed \end{theorem}

\begin{theorem} \label{Galois2}
Lifting Lemma: Given any objects 
$X \in \mathcal{C},\; Y \in \mathcal{C}$, and elements  
\mbox{$x \in FX$,}  \mbox{$y \in FY$,}  
if $lFix(x) \:\leq\: lFix(y)$ in $lAut(F)$, then there exist a unique 
arrow \mbox{$X \mr{f} Y$} in $\mathcal{C}$ such that $F(f)(x) = y$.

More generally, the following rule holds in $lAut(F)$
$$
\frac{\#[(X,\,\langle \!x_{0}\,|\,x_{1}\!\rangle )] \;\leq\; \#[(Y,\,\langle \!y_{0}\,|\,y_{1}\!\rangle )]}
  {\; \exists \; X \mr{f} Y \; F(f)(x_{0}) = y_{0}\,, \; F(f)(x_{1}) = y_{1}}
$$
\qed \end{theorem}
 
From \ref{lfiel}, \ref{transitivity2} and \ref{Galois2} follows by
easy categorical arguments:

\begin{theorem}[\textbf{fundamental theorem}] 
\label{teorema1}
Let $\mathcal{E}$ be a pointed connected atomic topos 
 $\cc{S} \mr{f} \mathcal{E}$, and $\mathcal{C} \subset
 \mathcal{E}$ be a pointed site $\mathcal{C} \mr{F}
 \cc{S}$ for $\mathcal{E}$, $\mathcal{E} = \mathcal{C}^{\sim}$,
$F = f^{*}|_{\cc{C}}$ (in the sense of \ref{atomicsite} above). Then, the
 lifting  $\mu F$ of $F$ (\ref{groupaction}) lands in the subcategory of
 transitive $G$-sets, 
$\mathcal{C} \mr{\mu F} t\mathcal{B}G$, for the localic group $G =
 lAut(F)$, and the induced morphism of
 topoi \mbox{$\cc{B}G \mr{\mu f} \cc{E}$}  is an equivalence.

\qed \end{theorem}

Actually, considering all the points (and with exactly the same proof that
theorem \ref{teorema0000}) this theorem yields:

\begin{theorem}[\textbf{fundamental theorem}] 
\label{Galois22}
Let $\mathcal{E}$ be any  pointed connected atomic topos. Then the
canonical geometric morphism \mbox{$\cc{B}(\cc{G}^{op}) \mr{} \cc{E}$}
is an equivalence, where $\cc{G}$ is the localic groupoid of
points $\cc{G} = l\cc{P}oints(\cc{E})$ defined in section \ref{lpoints}. 
\qed \end{theorem}


\section{Summary of the representation theorems and final conclusions} \label{final}

In this section we summarize an analysis of the results of this paper
and make some comments on Joyal-Tierney generalization of
Grothendieck's galois theory and its relation to the galois theory of
Galois topoi. We also treat the non connected theory, which follows
trivially from the connected case (as opposed to the groupoid formulation
from the group formulation within the connected theory).


\subsection{Comments on Joyal-Tierney galois theory}

In \cite{JT} Joyal-Tierney 
develop galois theory in a new way. Classifying topoi are explicitly
described as descent topoi. For them, the fundamental
theorem of galois theory states that \emph{(*) open
  surjections are geometric morphisms of effective descent}. 

The
fundamental theorems of previous galois theories follow because these
theorems are statements about a point, and this point is an open
surjection.  

It also follows an unpointed
theory of representation for a completely arbitrary topos in terms of
localic groupoids, theorem \mbox{Ch.VIII 3. Theorem 2.}, which
states that any topos is the classifying topos of a localic
groupoid. This theorem is dependent 
on change of base techniques. An important difference here with
previous galois theories is that the geometric morphism which is
proved to be of effective descent \emph{is not part or is not
  canonically associated to the data}. As a consequence, the groupoid
is not associated to the topos in a functorial way. We think that part of
Joyal-Tierney theory describes different phenomena from the one that
concerns the galois theory of Galois topoi, either pointed or unpointed. 

In \cite{JT} one also finds the representation theorem for pointed
atomic topoi \mbox{Ch.VIII 3. Theorem 1.} It is worth to
notice that Theorem 1 as such is not a particular
instance of Theorem 2. The proof of Theorem 2 goes in two steps. Step 1:
construct an open spatial (or localic) cover, which is the part that
does not corresponds to galois theory. Step 2:  using this cover
construct the localic groupoid that proves the statement by the
theorem (*) quoted above. In Theorem 1 
only the second step is used, the first one is already part of the
data (the given point is the cover), and as such it is canonical. The
atomic topoi with enough points have a canonical open spatial cover,
namely, the discrete localic space of all the points (only one is
necessary if the topos is connected), and it can be seen that the
construction in Step 2 yields the localic groupoid of points 
(defined in \ref{lpoints} above). The recipe given in section 3 of
chapter VII  for Step 1, applied to an atomic topos with enough points,
does not yield the discrete cover given by the points. 

We can say that  
pre-Joyal-Tierney galois theory cover Step 2 (and thus it suffices to
state and prove Theorem 1), as it has been shown in \cite{D}. While Step 1
(and thus Theorem 2) goes beyond.  


\subsection{The non connected theory}

Locally connected topoi are sums of connected locally connected
ones. Generally, because of this, it is enough to prove results for the
connected case. In \cite{G1} V 9. the non connected theory is left to
the reader (``Nous en laissons le detail au lecteur''). However, in
\cite{G3}, locally connected (but not connected) Galois topoi are
considered under the name \emph{``Topos Multigaloisiennes''}. There the
topos are supposed to have enough points.

\begin{definition}
A \emph{Multigalois topos} is a locally connected topos generated by
its Galois objects, or, equivalently, it is a sum of Galois topoi.

\end{definition}
Let $\cc{E}$ be any locally connected topos. The following theorems follow
by decomposing $\cc{E}$ as a sum of connected topoi, and proving the
statements for connected topoi (which we indicate in between
parenthesis). The implication chain is best
understood if performed by the increasing cycle permutation.
 
\vspace{1ex}

From the results in section \ref{Galois' galois theory}  we have:  
\begin{theorem} [\textbf{discrete case}] \label{discrete case}
The following are equivalent:

\indent 2. $\cc{E}$ is a locally simply-connected (Galois) Multigalois topos.
             
\indent 3. $\cc{E}$ is (connected) atomic with enough 
essential points.

\indent 4. The canonical geometric morphism $\cc{B}(\cc{G}) \mr{\mu}
\cc{E}$ is an equivalence, where $\cc{G}$ is the (connected) discrete
groupoid of points $\cc{G} = \cc{P}oints(\cc{E})$.

\indent 5. $\cc{E}$ is the classifying topos of any (connected) discrete
groupoid.
\qed \end{theorem} 

From the results in section \ref{Grothendieck's galois theory} (and
\ref{GLCisgalois2}) we have:
\begin{theorem}[\textbf{pointed prodiscrete case}] \label{prodiscrete case}
The following are equivalent:

\indent 1. $\cc{E}$ has enough points and it is (connected) generated by its
locally constant objects. 

\indent 2. $\cc{E}$ has enough points and it is a (Galois) Multigalois
topos.

\indent 3. $\cc{E}$ is (connected) atomic with enough proessential
points.

\indent 4. The canonical geometric morphism $\cc{B}(\cc{G}) \mr{\mu}
\cc{E}$ is an equivalence, where $\cc{G}$ is the (connected) localic groupoid
of points $\cc{G} = l\cc{P}oints(\cc{E})$, which in
this case has prodiscrete ``hom-spaces''.

\indent 5. $\cc{E}$ is the classifying topos of any (connected)
localic groupoid with discrete space of objects and prodiscrete
``hom-spaces''.
\qed \end{theorem} 
 
From the results in section \ref{unpointedsection} (and
\ref{GLCisgalois2}) we have:
\begin{theorem}[\textbf{unpointed prodiscrete case}] \label{unpointed
prodiscrete case}
The following are equivalent:

\indent 1. $\cc{E}$ is (connected) generated by its locally constant objects. 

\indent 2. $\cc{E}$ is a (Galois) Multigalois topos.

\indent 4. The canonical geometric morphism $\cc{B}(\cc{G}) \mr{\mu}
\cc{E}$ is an equivalence, where $\cc{G}$ is the (connected) localic groupoid
of phantom points $\cc{G} = ph\cc{P}oints(\cc{E})$, which is
prodiscrete (by definition).

\indent 5. $\cc{E}$ is the classifying topos of any (connected)
prodiscrete localic groupoid.
\qed \end{theorem} 

From the results in section \ref{localicgalois} we have:
\begin{theorem}[\textbf{pointed localic case}] \label{localic case}
The following are equivalent:

\indent 3. $\cc{E}$ is (connected) atomic with enough points.

\indent 4. The canonical geometric morphism $\cc{B}(\cc{G}) \mr{\mu}
\cc{E}$ is an equivalence, where $\cc{G}$ is the (connected) localic groupoid
of points $\cc{G} = l\cc{P}oints(\cc{E})$.

\indent 5. $\cc{E}$ is the classifying topos of any (connected)
localic groupoid with discrete space of object.
\qed \end{theorem} 

We see that the theorem that should be labeled \emph{unpointed localic
case} is missing. It should concern an
arbitrary, (even connected but may be pointless) atomic topos. Of
course \mbox{\cite{JT} Ch.VIII 3. Theorem 2}  applies to such a topos, but
it is a far too general theorem. The localic groupoid $\cc{G}$ 
such that $\cc{E} \cong \cc{B}\cc{G}$ is not identified
and can not be
considered to be (as far as we can imagine) the groupoid of points of
the topos.  In \cite{M3} Moerdiejk investigates Joyal-Tierney theorems and
establishes that atomic topoi are characterized by the fact that the
localic groupoid can be so chosen that the map 
$G_1 \mr{(d_0,\, d_1)} G_0 \times G_0$ is open (\cite{M3} 4.7 c). However we
are still far
from any canonicity for the groupoid.  We need a theorem (still
unknown) which, in particular, in the
presence of a point, yields \cite{JT} Ch.VIII 3. Theorem 1. We still
not know how to define the groupoid of phantom points for a general
atomic topos (as we do for a general Galois topos). The solution is
not given by Joyal-Tierney's generalization of Grothendieck's galois
theory, and it is not its purpose either. 
An \emph{unpointed localic galois theory} (compressing the unpointed
prodiscrete galois theory as well as the pointed localic galois
theory) is yet to be developed.




\section{Appendix. Galois theory of covering topoi} \label{locallyconstantsection}

We shall denote by $\gamma$ the structure morphism of any topos $\cc{E}$,
$\cc{E} \mr{\gamma} \cc{S}$. A topos $\cc{E}$ is said to be \emph{locally
connected} if the inverse image functor $\gamma^{*}$ is essential, that
is, if it has itself a left
adjoint denoted $\gamma_{\,!}$ (the set of connected components). A
topos $\cc{E}$ is said to be \emph{connected} if the inverse image
functor $\gamma^{*}$ is full and faithful. If $\cc{E}$ is connected
and locally connected clearly 
$\gamma_{\,!}\gamma^{*} = id$. The
reference for connected and locally connected topoi is 
\cite{G2}, Expos\'e IV, 4.3.5, 4.7.4, 7.6 and 8.7. A geometric
morphism  $\cc{E} \to \cc{F}$ is said to be a \emph{locally connected
morphism}  if the
topos $\cc{E}$ considered as an $\cc{F}$-\emph{Topos} is locally
connected. This relative version was introduced in    
 \cite{BP} under the name $\cc{F}$-\emph{essential}, see also 
the appendix of \cite{M1}. Recall that a \emph{connected atomic
topos} is a connected, locally connected and boolean topos.  
For atomic topoi and atomic sites the reference is \cite{BD}, 
also see \cite{JT}.

\vspace{1ex}

A \emph{covering} of a topos $\cc{E}$ is a geometric morphism of the form 
$\cc{E}/X  \mr{} \cc{E}$, with $X$ a \emph{locally constant object}. 

\vspace{2ex}

\subsection{Locally constant objects}
\indent \vspace{1ex}  \label{LCOsection}

We recall now the definition of locally constant object in an
arbitrary topos given in SGA4, Expos\'e IX (see also \cite{BL} 
where this definition is considered over an arbitrary base topos). 

\begin{definition} \label{lcc} An object $X$ of
a topos $\cc{E} \mr{\gamma} \cc{S}$ is said to be $U$-split, for a
cover $U = \{U_{i}\}_{i \in I}$ (i.e. epimorphic family $U_{i} \to 1$), if it
becomes constant on each $U_{i}$. That is, if there exists family of sets 
$\{S_{i}\}_{i \in I}$ and isomorphisms in $\cc{E}$,
\mbox{$\{\gamma^{*}S_{i} \times U_{i} \mr{\theta_{i}}  X \times
U_{i}\}_{i \in I}$} over $U_i$.
 We say that $X$ is \emph{locally
constant} if $X$ is $U$-split for \emph{some} cover U in $\cc{E}$.
\end{definition}

It is often convenient to identify a family with 
a function $\beta\colon
S \to I$, $S = \sum_{i} S_{i}$. We abuse the language and write also
$U$ for the coproduct $U = \sum_{i} U_{i}$. Notice that there is a map
$\zeta \colon U \to {\gamma}^{*}I = \sum_{i} 1$, and in this way a 
cover as above is given by such a map with $U \to 1$ epimorphic. In this
notation the family
of isomorphisms $\theta_{i}$ is the same thing as an isomorphism
$\theta \colon {\gamma}^{*}S \times_{{\gamma}^{*}I} U \to X \times U$
over $U$.

\vspace{1ex}
  
When the topos is connected a
classical (in the theory of topological coverings spaces) connectivity
argument shows that all the sets
$S_{i}$ can be considered equal (see \cite{BDLSCT} for a proof of this
in the topos context). If the topos is not connected this ``single set''
concept is clearly not equivalent (we can have a different set for
each connected component of 1), and not the right one, as it has been
observed in \cite{BARR}. We comment to the reader
interested in the relative theory over an arbitrary base topos that the
connectivity argument depends on the excluded middle. Based on this,
even when the topos is connected, in the
relative case the ``single set version'' of the notion will not be
equivalent to the ``family version'' of SGA4 or \cite{BL}, even for 
connected topoi (unless the base topos is boolean). I thank here A. Kock 
for some fruitful conversations on the connectivity argument in 
\cite{BDLSCT}, which led me to conjecture that if the argument is valid 
for any connected locally connected topos $\cc{E} \to \cc{S}$, then 
necessarily $\cc{S}$ has to be boolean.

\vspace{1ex}

 (*) Assume now that $\cc{E}$ is locally connected. In this
case, by enlarging $I$, we can always consider $I = {\gamma}_{!}U$. In fact, 
the connected components of \mbox{$U = \sum_{i} U_{i}$} are the
connected components 
of the $U_{i}^{\;'s}$. Repeat then the same set $S_i$ for each connected 
component of $U_i$. The map 
$\zeta \colon U \to {\gamma}^{*}{\gamma}_{!}U$ results the unit of
the adjunction ${\gamma}_{!} \dashv {\gamma}^{*}$ at $U$.

\vspace{1ex}


We consider now a locally connected topos. Given a cover $U$, the full subcategory $Split(U)$ of objects
split by $U$ is a topos (see \cite{G2}, \cite{BDLSCT}), in fact a
quotient topos with inverse image given by the inclusion. Obviously
the adjunction ${\gamma}_{!} \dashv {\gamma}^{*}$ restricts to
$Split(U)$, and so this topos is locally connected. It is immediate to
check that it is boolean (given $Z \hookrightarrow X$ in $Split(U)$,
then $Z \times U \,\hookrightarrow \, X \times U$ has a complement and
this implies that $Z \hookrightarrow X$ has one, (see \cite{BDLSCT}). 
Thus $Split(U)$ is an
atomic topos, and clearly, if $\cc{E}$ is
connected, so is $Split(U)$. These results are derived by elegant and
simpler (but indirect, and less convincing in a way) arguments in \cite{BM}.  


Given any topos, the covers $U \to 1$ epi form a (co)
filtered poset $Cov(\cc{E})$ if we define 
$U \leq V \;\iff\; \exists \; U \to V$. This
poset has a small cofinal subset. In fact, as observed in \cite{BD},
the irredundant sums of generators with global support, of which there
is only a set, are cofinal.

If $U \leq V$, it is immediate to check that $Split(V) \subset
Split(U)$, and that the inclusion is the inverse image of a geometric
morphism of topoi. In this way we have a filtered inverse limit of
topoi. Clearly the inverse limit site for this topos is the full
subcategory  of all connected locally constant objects of
$\cc{E}$, and the topos, that we denote $GLC(\cc{E})$, as a category,
is the full subcategory of objects generated by the locally constant objects.
It follows that the inclusion is the inverse image of a geometric
morphism, and $\cc{E} \to GLC(\cc{E})$ is a quotient topos. Again, the
adjunction ${\gamma}_{!} \dashv {\gamma}^{*}$ obviously
restricts to $GLC(\cc{E})$, and so this topos is 
locally connected. In the same way that for $Split(U)$ it is immediate
to check that it is boolean, thus it is an atomic topos, and if
$\cc{E}$ is connected, so is $GLC(\cc{E})$. We resume
this situation in the following diagram:
$$Split(U)  \ml{} Split(V)  \;\;\;\;\; \cdots  \ml{} GLC(\cc{E})
                                               \ml{} \cc{E}$$
where $GLC(\cc{E})$ is a filtered  inverse limit of topoi indexed by
the poset $Cov(\cc{E})$.


\vspace{1ex}

In \cite{B} a push-out topos is considered in order to define
categories of locally constant objects.
$$
\xymatrix
        {
         \cc{E}/U \ar[r]^{\varphi_U} \ar[d]^{\rho_U} 
        & \cc{E} \ar[d]^{\sigma_{U}}
        \\
         \cc{S}/\gamma_{!}U \ar[r]^{f_U} 
        & \cc{P}_U
        }
$$
where $\rho_{U}$ and $\varphi_U$ are given by 
$\rho_{U}^{*}(S \to \gamma_{!}U) \;=\; {\gamma}^{*}S \;
{\times}_{{\gamma}^{*}{\gamma}_{!}U} \; {U}$  and
\mbox{$\varphi_{U}^{*}(X) \;=\; X \;{\times} \; {U}$}. It is well
known that the geometric morphism $\varphi_U$ is locally connected,
and it can be
checked that the geometric morphism $\rho_{U}$ is connected and
locally connected. It follows that $f_U$ is locally connected  and
that $\sigma_{U}$ is connected locally connected (see \cite{B} lemma 2.3).

Consider the construction of push outs of
topoi. An object of $\cc{P}_{U}$ is a $3$-tuple $\langle X, S \to \gamma_{!}U,
\theta\rangle $, with $\theta\colon{X \times U}\to {\gamma^{*}}S
\times_{\gamma^{*}\gamma_{!}U} U$  an isomorphism over $U$, and a morphism
 \mbox{$\langle X, S \to \gamma_{!}U, \theta\rangle  \;\to\;  \langle X', S' \to \gamma_{!}U, 
\theta'\rangle $}
is determined by a pair of morphisms $f \colon X \to X'$ and
$\alpha \colon S \to S'$,
the latter over $\gamma_{!}U$, compatible with the isomorphisms $\theta$ and
$\theta'$. The functor $\sigma_{U}^{*}$ is the projection
functor $\cc{P}_{U} \to \cc{E}$, which is then fully faithful (the
reader can also check by direct inspection that $f \colon X \to X'$
determines $\alpha \colon S \to S'$ once we assume that $X$ and $X'$
are part of the data of $U$-locally constant objects in $\cc{E}$). 

By considerations made above (*), the essential image of
$\sigma_{U}^{*}$ is the full subcategory $Split(U)$, thus we have:
   
\begin{proposition} \label{GG_U}
Given any locally connected topos $\cc{E}$ and a cover $U$,
 the push-out topos
$\cc{P}_{U}$ is equivalent (as a category) via the full and faithful
projection functor \mbox{$\sigma_{U}^{*}:\;
\cc{P}_{U}  \mr{\cong} Split(U) \subset  \cc{E}$} to the full
subcategory $Split(U)$. Thus $\cc{P}_{U}$ and $Split(U)$ are
equivalent topoi.
\qed \end{proposition}

The morphism $f_U$ 
is actually a family of points indexed by
$\gamma_{!}U$, and since it is locally connected, all these points are
essential. The inverse image of $f_U$ is given by 
$f_{U}^{*}\langle X, S \to \gamma_{!}U, \theta\rangle \; =\;  S \to \gamma_{!}U$. We can 
prove directly:

\begin{proposition} \label{GG_Upoints0}
Given any locally connected topos $\cc{E}$ and a cover $U$, for each 
$i \in \gamma_{!}U$ the composite, denoted $f_i$, of the
corresponding point of $\cc{S}/\gamma_{!}U$ with $f_U$ is an
essential point of $\cc{P}_{U}$. 
\end{proposition}
\begin{proof}
The inverse image of $f_i$ is given by $f_{i}^{*}\langle \!X, S \to \gamma_{!}U,
\theta\!\rangle  \;=\; S_i$. It follows by the construction of inverse limits in the
push-out topos that $f_{i}^{*}$ preserves all inverse limits, that is,
it is essential. 
\end{proof}

Given any connected locally connected topos $\cc{E}$, it follows then from
proposition \ref{galois1} that
\emph{all the results of section \ref{Galois' galois theory} 
apply to the topoi $\cc{P}_U$ and $Split(U)$}. 
From proposition \ref{Asplit} we have the following
important fact (existence of Galois closure):

\begin{proposition} \label{U=A}
Given any connected locally connected topos $\cc{E}$ and any cover
$U$, $Split(U) = Split(A)$, and $\cc{P}_U \cong \cc{P}_A$, for $A$ any 
representing object (necessarily a Galois object (\ref{GOsection}), thus
in particular 
a connected cover) of one of the points $f_i$.
\qed \end{proposition}

\begin{proposition} \label{GG_Upoints1}
In the situation of proposition \ref{GG_Upoints0}, any point $g$ is
isomorphic to some $f_i$.  
\end{proposition}
\begin{proof}
Let $g$ be any 
point. Since the family $U_i \to 1$ is epimorphic it follows that
there is (at least) one $i$ with $g^{*}U_i \to  1$ epimorphic, thus 
$g^{*}U_i \not= \emptyset$. Given any
object $X$ split by
$U$, since $g^{*}\gamma^{*} = id$ we have an isomorphism
$S_{i} \times g^{*}U_{i} \mr{\cong}  g^{*}X \times
g^{*}U_{i}$ over $g^{*}U_i$. This clearly implies $g^{*}X \cong S_i$. 
\end{proof}

It follows then from proposition \ref{basica} that when the topos is
connected, all the points $f_i$ are isomorphic.

Furthermore, given any locally connected topos $\cc{E}$ and a connected cover
$U$ (notice that this forces $\cc{E}$  to be a
connected topos) we have:

\begin{proposition} \label{GG_Upointsunique}
Given any connected locally connected topos $\cc{E}$ and a connected 
cover $U$, the topos $\cc{P}_{U}$ has a canonical essential point 
$f = f_U$, and any other point is isomorphic to $f$.
\qed \end{proposition} 

Here it is important to stress the fact that although there is a
canonical (geometrical morphism) equivalence 
$Split(U) \mr{\cong} \cc{P}_{U}$, the topos
$Split(U)$  \emph{does not have a canonical point} since the
equivalence does not have a canonical inverse.  

\vspace{2ex}

\subsection{Galois objects} \indent \vspace{1ex} \label{GOsection}

Recall that a non-empty connected object $A$ (called a molecule
in \cite{BP}) in a topos $\cc{E}$
is said to be a  
\emph{Galois object} if it is an $Aut(A)$-torsor. That is, if $A \to
1$ is an epimorphism, and the 
canonical morphism \mbox{$A \times \gamma^{*} Aut(A) \mr{}  A \times A$} is an
isomorphism. Clearly, any Galois object is a locally constant
object. Notice also that non-empty connected objects in atomic topoi
are \emph{atoms} in the sense that the lattice of subobjects is $2$.   

After Grothendieck's ``Categories Galoisiennes'' of \cite{G1}
and Moerdiejk ``Galois Topos'' of \cite{M2}, we state the following
definition:

\begin{definition} \label{galoistopos}
A Galois Topos is a connected locally connected topos generated by its
Galois objects. 
\end{definition}
 
Notice that unlike \cite{G1} and \cite{M2} we do not require the topos
to be pointed. Although all the applications concern pointed topoi, it
is still interesting to notice that the basic theory can be
developed without this assumption, as it has been shown in \cite{B},
\cite{BM} and  \cite{K}.

Since Galois objects are locally constant, clearly
the canonical morphism gives an equality of topoi $\cc{E} =
GLC(\cc{E})$. In particular, Galois topoi are atomic.

\vspace{1ex}

From the equation $\gamma^* f^* = \gamma^*$ and the fact that inverse
image $f^*$ of a connected geometric
morphism \mbox{$\cc{E} \mr{f} \cc{F}$} is a full and faithful left-exact
functor, it immediately follows that an object $A$ in $\cc{F}$ is an
$Aut(A)$-torsor if and only if the object $f^{*}A$ in $\cc{E}$ is an
$Aut(f^{*}A)$-torsor. Furthermore, using now in addition that $f^*$
preserves (in particular)  binary coproducts, it easily follows that 
$A$ is a non empty connected object if and only if $f^{*}A$ is so. In
conclusion we have:

\begin{proposition}
Given any connected geometric morphism $\cc{E} \mr{f} \cc{F}$, an
object $A \in \cc{F}$ is a Galois object if and only if $f^{*}A \in
\cc{E}$ is so. 
\qed\end{proposition}   

\begin{proposition} \label{inverselimitgalois}
Any filtered inverse limit of Galois topoi and connected locally
connected geometrical morphisms is a Galois topos.
\end{proposition}
\begin{proof}
That the inverse limit topos
is connected and locally connected is proved in \cite{M1}. Consider
now the corresponding colimit of sites as in section
\ref{inverselimitsection}. By construction of the inverse limit
site and the previous proposition it follows that the
Galois objects generate the inverse limit topos. 
\end{proof}

From this proposition and proposition \ref{U=A} it follows:

\begin{theorem} \label{GLCisgalois}
Given any connected locally connected topos $\cc{E}$, the topos
$GLC(\cc{E})$ is a Galois topos. 
\qed \end{theorem}

Notice  that it follows the equality 
$GLC(GLC(\cc{E})) = GLC(\cc{E})$ (any locally constant object is split
by a locally constant cover), fact which is not evident by definition.

\begin{corollary} \label{GLCisgalois2}
Given any connected locally connected topos $\cc{E}$, then $\cc{E}$ is
a Galois topos if and only if $\cc{E}$ is generated by its locally
constant objects if and only if $\cc{E} = GLC(\cc{E})$.
\qed \end{corollary}

Let $GCov(\cc{E})$ be the subposet of $Cov(\cc{E})$ whose objects are
Galois objects (necessarily covers). Although \emph{it is not cofinal}, it is also 
filtered. In fact, given two Galois objects $A$, $B$, consider the
Galois object $C$ such that $Split(A \times B) = Split(C)$ given by
proposition \ref{U=A}. We have:
 
\begin{proposition} \label{filtergalois}
Any Galois topos is the filtered inverse limit of the topoi
$Split(A)$, $A \in GCov(\cc{E})$. The inverse limit site is the
filtered union of the full subcategories $cSplit(A) \subset \cc{E}$ of
connected objects split by $A$.
\qed \end{proposition}

Given any connected locally connected topos we can now synthesize the
situation in the following diagram:
$$Split(A)  \ml{} Split(B)  \;\;\;\;\; \cdots  \ml{} GLC(\cc{E})
                                               \ml{} \cc{E}$$


where $GLC(\cc{E})$ is a filtered  inverse limit of topoi indexed by the
poset $GCov({\cc{E}})$ whose objects are the Galois covers. $\cc{E}$ is
a Galois topos if and only if the left-most arrow is the equality. 

\vspace{1ex}

When $\cc{E}$ is a pointed topos, clearly $GLC(\cc{E})$ is also
pointed, thus
it follows that \emph{all the results of section  \ref{Grothendieck's
galois theory} apply to the topos $GLC(\cc{E})$}.

\vspace{1ex}

The original definition of Galois object given in \cite{G1} was
relative to a point of the topos. However, that point was surjective,
and it is easy to check:
 
\begin{proposition} \label{chargalois}
Let $\mathcal{E}$ be a topos furnished with a surjective point (meaning 
the inverse image functor reflects isomorphisms),
 $\cc{S} \mr{f} \mathcal{E}$. Then, a non-empty connected
object $A$ is a Galois object if and only if there exists $a \in f^{*}A$
so that the map $Aut(A) \mr{a^{*}} f^{*}A$, defined by $a^{*}(h) =
f^{*}(h)(a)$ is a bijection (the same holds then for any other element $b
\in f^{*}A$). 
\qed \end{proposition}

Notice that this characterization of Galois objects is word by word
equal to the definition of normal extension in the classical Artin's
interpretation of galois theory.



\end{document}